\documentclass[12pt,reqno]{amsart}
\usepackage{amsmath}
\usepackage{amsfonts}
\usepackage{amssymb}
\usepackage{caption}
\usepackage{latexsym}
\usepackage{color}
\usepackage{hyperref}
\usepackage{url}
\usepackage{enumerate}
\usepackage{arydshln} 
\usepackage{multicol} 

\renewcommand{\proof}{\noindent{\it Proof.\ \ }}
\renewcommand{\qed}{\ifmmode\square\else\nolinebreak\hfill
$\Box$\fi\par\vskip12pt}

\renewcommand\a{\alpha}    

\newcommand\Ga{\mathrm{\Gamma}}   \newcommand\Ome{{\it \Omega}}
\newcommand\Sig{{\it \Sigma}}

\newcommand\A{\mathrm{A}}  \newcommand\B{\mathrm{B}}  \newcommand\C{\mathbf{C}}   \newcommand\D{\mathrm{D}}
  \newcommand\F{\mathrm{F}}  \newcommand\G{\mathrm{G}}   
\newcommand\K{\mathsf{K}}  \newcommand\M{\mathrm{M}}     
\renewcommand\P{\mathrm{P}}\newcommand\Q{\mathrm{Q}}     

\newcommand\ZZ{\mathbb{Z}}     \newcommand\BB{\mathcal{B}}        \newcommand\PP{\mathcal{P}}

\newcommand\mz{\mathbb{Z}}

          \newcommand\Aut{\mathrm{Aut}}
   \newcommand\Cay{\mathrm{Cay}}      
      \newcommand\Inn{{\mathrm Inn}}         
\newcommand\Out{\mathrm{Out}}    \newcommand\PG{\mathrm{PG}}        
\newcommand\soc{\mathrm{soc}}    \newcommand\Sy{\mathrm{S}}         
          \newcommand\Mult{\mathrm{Mult}}
\newcommand\diam{\mathrm{diam}}            \newcommand\la{\langle}
\newcommand\ra{\rangle}

\newcommand\AGL{\mathrm{AGL}}                  
       \newcommand\AGammaL{\mathrm{A\Gamma L}}

\newcommand\GammaL{\mathrm{\Gamma L}}          
\newcommand\GL{\mathrm{GL}}                    
                    
\newcommand\PGL{\mathrm{PGL}}                  \newcommand\PGammaL{\mathrm{P\Gamma L}}
                  
\newcommand\PSL{\mathrm{PSL}}

          \newcommand\SL{\mathrm{SL}}

\newtheorem{theorem}{Theorem}[section]%
\newtheorem{lemma}[theorem]{Lemma}%
\newtheorem{proposition}[theorem]{Proposition}%
\newtheorem*{conjecture}{Conjecture}%

\begin{document}

\title[The classification of two-distance transitive dihedrants]
{The classification of two-distance transitive dihedrants}

\thanks{$^*$Corresponding author. }
\thanks{1991 MR Subject Classification 05C25, 20B15, 20B30.}
\thanks{This work was partially supported by the National Natural Science Foundation of China (12331013, 12311530692, 12271024, 12161141005, 12071023, 12301461) and the
111 Project of China (B16002).}

\author[J.-J. Huang, Y.-Q. Feng, J.-X. Zhou, F.-G. Yin]{Jun-Jie Huang, Yan-Quan Feng$^*$, Jin-Xin Zhou, Fu-Gang Yin}
\address{Jun-Jie Huang, Yan-Quan Feng, Jin-Xin Zhou\\
School of mathematics and statistics\\
Beijing Jiaotong University\\
Beijing \\
100044, P. R. China}
\address{Fu-Gang Yin\\
School of Mathematics and Statistics\\
Central South University\\
Changsha, Hunan \\
410083, P.R. China}

{\email{20118006@bjtu.edu.cn(J.-J. Huang), yqfeng@bjtu.edu.cn (Y.-Q. Feng), \linebreak jxzhou@bjtu.edu.cn (J.-X. Zhou), 18118010@bjtu.edu.cn(F.-G. Yin)}
\maketitle

\begin{abstract}
  A vertex transitive graph $\Gamma$ is said to be {\em $2$-distance transitive} if for each vertex $u$, the group of automorphisms of $\Gamma$ fixing the vertex $u$ acts transitively on the set of vertices at distance $1$ and $2$ from $u$, while $\Gamma$ is said to be {\em $2$-arc transitive} if  its automorphism group is transitive on the set of $2$-arcs. Then $2$-arc transitive graphs are $2$-distance transitive.  The classification of $2$-arc transitive Cayley graphs on dihedral groups was given by Du, Malni\v{c} and Maru\v{s}i\v{c} in [Classification of 2-arc-transitive dihedrants, J. Combin. Theory Ser. B 98 (2008), 1349--1372]. In this paper, it is shown that a connected 2-distance transitive Cayley graph on the dihedral group of order $2n$ is either $2$-arc transitive, or isomorphic to the complete multipartite graph $\K_{m[b]}$ for some $m\geq3$ and $b\geq2$ with $mb=2n$.
\end{abstract}

\qquad {\textsc k}{\scriptsize \textsc {eywords.} 2-distance transitive, Cayley graph, dihedral group, quasiprimitive group} {\footnotesize}

\section{Introduction}

In this paper, all graphs are finite, simple, connected and undirected.
For a graph $\Ga$, we use $V(\Ga)$, $E(\Ga)$, $A(\Ga)$ and $\Aut(\Ga)$ to denote its vertex set, edge set, arc set and automorphism group, respectively.
If $\Aut(\Ga)$ is transitive on $V(\Ga)$, $E(\Ga)$ or $A(\Ga)$, then $\Ga$ is said to be {\em vertex transitive}, {\em edge transitive} or {\em arc transitive}, respectively.
Let $u$ and $v$ be two distinct vertices of $\Ga$.
The {\em distance} between $u$ and $v$ in $\Ga$ is the smallest length of paths between $u$ and $v$, denoted by $d_\Ga(u,v)$, while the {\em diameter} of $\Ga$ is the maximal value of the distance of all vertex pairs, denoted by $\diam(\Ga)$.
For a positive integer $i$, denote by $\Ga_i(u)$ the set of vertices at distance $i$ with vertex $u$ in $\Ga$.
In particular, $\Ga_1(u)$ is simply denoted by $\Ga(u)$.
Let $\Ga$ be a graph admitting a vertex transitive group $G$ of  automorphisms of $\Ga$.
Then $\Ga$ is said to be {\em $(G,s)$-distance transitive}, if for each vertex $u\in V(\Ga)$, the stabilizer $G_u$ is transitive on $\Ga_i(u)$ for all $i\leq s$.
If $s=\diam(\Ga)$, then $\Ga$ is said to be {\em $G$-distance transitive}.
When $G=\Aut(\Ga)$, a $(G,s)$-distance transitive graph or $G$-distance transitive graph will be simply called an {\em $s$-distance transitive graph} or a {\em distance transitive graph}, respectively.

Let $G$ be a finite group and let $S$ be a subset of $G\setminus\{1\}$ such that $S=S^{-1}$ and $G=\la S\ra$. The {\em Cayley graph} on $G$ with respect to $S$ is defined to be the graph with vertex set $G$ and edge set $\{\{g,sg\}\mid g\in G,s\in S\}$, denoted by $\Cay(G,S)$.
It is widely known that a graph $\Ga$ is a Cayley graph on $G$ if and only if $\Aut(\Ga)$ has a regular subgroup isomorphic to $G$, see~\cite[Lemma 16.3]{Biggs}. Given $g\in G$, the right multiplication $R(g): x\mapsto xg$, $x\in G$, is an automorphism of $\Cay(G,S)$, and $R(G):=\{R(g)\ | \ g\in G\}$ is a regular group of automorphisms of $\Cay(G,S)$, that is, the image of $G$ under its right regular representation. A Cayley graph $\Cay(G,S)$ is said to be {\em normal} if $R(G)$ is a normal subgroup of $\Aut(\Cay(G,S))$. Furthermore, a Cayley graph $\Cay(G,S)$ is called a {\em dihedrant} if $G$ is a dihedral group.

The investigation of distance-transitive graphs is one in which the theoretical developments in algebraic graph theory can be traced back to the 1980s. Many researchers have contributed significantly to this field and have achieved important research results. For a comprehensive overview, refer to the book~\cite{BCN} and the references therein.
Notably, all distance-transitive graphs are $2$-distance transitive, but the converse is not true.
In fact, an infinite family of normal Cayley graphs on the extraspecial $p$-group of order $p^3$ and exponent $p$ with $p$ an odd prime was constructed in~\cite{HFZ}. These graphs are $2$-distance transitive but are neither distance transitive nor $2$-arc transitive.
The extensive study of $2$-distance transitive graphs has gained momentum in recent years, as evidenced by works such as those cited in~\cite{CJL,CJS,HFZ,JT,JWZ}.
This paper aims to contribute to this evolving topic.

Let $u,v$ and $w$ be three distinct vertices of a graph $\Ga$. The triple $(u,v,w)$ is called a {\em $2$-arc} if $v$ is adjacent to both $u$ and $w$, and $\Ga$ is said to be {\em $2$-arc transitive} if $\Aut(\Ga)$ is transitive on $V(\Ga)$ and on the set of $2$-arcs of $\Ga$.
It is important to note that all $2$-arc transitive graphs are $2$-distance transitive, but the reverse is not necessarily true.
For example, the complete multipartite graphs $\K_{m[b]}$ with $m\geq 3$ parts, each containing $b\geq 2$ vertices, serve as instances where the graphs are $2$-distance transitive but not $2$-arc transitive.
There has been a significant body of work on the $2$-arc transitivity of graphs in the literature, see \cite{ACX, DMM, LSS, Mar03, Mar06} for some examples.
When studying $2$-distance transitive graphs, particular interest lies in graphs that are not $2$-arc transitive.

The symmetry of Cayley graphs has attracted significant attention from researchers, with one noteworthy contribution being the completion of the classification of $2$-distance transitive Cayley graphs on cyclic groups. This work is documented in \cite{ACX,CJL}.
Subsequently, a natural problem arises: the classification of $2$-distance transitive Cayley graphs on dihedral groups. There have been fascinating classification results on dihedrants under special symmetries. For instance, the classification of $2$-arc transitive dihedrants was presented across a series of papers \cite{DMM, Mar03, Mar06}, and the classification of distance transitive or locally primitive dihedrants was detailed in \cite{MP} and \cite{Pan}, respectively.
Furthermore, the classifications of arc transitive or edge transitive Cayley graphs on the dihedral group $\D_{2p^n}$ where $p$ is a prime, were conducted in \cite{Kov} and \cite{SLZ}, respectively. The classifications of quasiprimitive or biquasiprimitive edge-transitive dihedrants were also provided in \cite{PYZH}.
Recently, Jin and Tan explored $2$-distance transitive dihedrants of girth $3$ in their work \cite{JT} and formulated the following conjecture (refer to \cite[Conjecture 1.4]{JT}):

\begin{conjecture}\label{conj}
A connected $2$-distance transitive dihedrant either is known $2$-arc transitive dihedrant, or is isomorphic to $\K_{m[b]}$ for some $m\geq3$ and $b\geq 2$, or $\G(2,p,(p-1)/4)$ for a prime $p$ with  $p\equiv1\pmod{8}$.
\end{conjecture}

In this conjecture, $\G(2,p,(p-1)/4)$ represents a family of arc transitive graphs of order $2p$ that were constructed in \cite[P. 199]{CO}. However, it can be demonstrated that $\G(2,p,(p-1)/4)$ is not $2$-distance transitive. This fact can be verified by examining certain graphs with small values of $p$, such as $p=17$ or $p=73$, which can be checked using computational tools like Magma \cite{Magma}.
In the paper, the conjecture mentioned is confirmed by excluding the inclusion of $\G(2,p,(p-1)/4)$.

\begin{theorem}\label{Thm-1}
Let $\Ga$ be a connected Cayley graph on the dihedral group $\D_{2n}$ of order $2n$ with $n\geq2$. Then $\Ga$ is $2$-distance transitive if and only if either $\Ga$ is a $2$-arc transitive graph, or $\Ga\cong\K_{m[b]}$ for some integers $m\geq3$ and $b\geq2$ with $mb=2n$.
\end{theorem}

\noindent{\bf Remark.}\ (1) By Theorem~\ref{Thm-1}, the classification of $2$-distance transitive dihedrants depends on the classification of $2$-arc transitive dihedrants. The study of $2$-arc transitive dihedrants was initiated by Maru\v si\v c in \cite{Mar03}, where a structural reduction theorem for $2$-arc-transitive dihedrants was established. Subsequently, the classification of $2$-arc transitive dihedrants was completed in \cite{DMM}, which was supplemented in \cite{Q-Du-Koolen} by the construction of new $2$-arc transitive dihedrants that were not covered in \cite{DMM}.

(2) A $2$-arc $(u,v,w)$ of a graph $\Ga$ is termed a {\em $2$-geodesic} if $u$ and $w$ are not adjacent, and an arc transitive graph $\Ga$ is said to be {\em $2$-geodesic transitive} if $\Aut(\Ga)$ is transitive on the set of $2$-geodesics. While all $2$-geodesic transitive graphs are $2$-distance transitive, the reverse implication does not hold true. This distinction is illustrated, for instance, in \cite[Theorem 1.2]{JDLP}.
As a consequence of Theorem~\ref{Thm-1}, a classification of connected $2$-geodesic transitive dihedrants can be derived, as outlined in \cite[Theorem 1.5]{JT}.

\section{Preliminaries}\label{Preliminaries}

In this section, we provide definitions and results on groups and graphs that will be utilized in our analysis. For a positive integer $n$, denote by $\ZZ_n$ the additive group of integers modulo $n$, by $\ZZ_n^*$ the multiplicative group of numbers coprime to $n$ in $\ZZ_n$, and by $\D_{2n}$ the dihedral group of order $2n$. For a prime $p$ and a positive integer $r$, denote by $\ZZ_p^r$ the elementary abelian group of order $p^r$. For two groups $A$ and $B$, denote by $A\times B$ the direct product of $A$ and $B$, by $A:B$ a semidirect product of $A$ by $B$, and by $A.B$ an extension of $A$ by $B$. For a group $G$, denote by $Z(G)$ the center of $G$, and by $\soc(G)$ the socle of $G$, that is, the product of all minimal normal subgroups of $G$. Denote by $\C_n$ the cycle graph of length $n$, $\K_n$ the complete graph of order $n$, $\K_{n,n}$ the complete bipartite graph of order $2n$, $\K_{n,n}-n\K_2$ the subgraph of $\K_{n,n}$ minus a matching, and $\K_{m[b]}$ the complete multipartite graph consisting $m\geq 3$ parts of size $b\geq 2$.
\subsection{Group theoretical results}

An extension $G=N. H$ of $N$ by $H$ is termed a {\em central extension} if $N\leq Z(G)$, and a central extension $G=N. H$ is referred to as a {\em covering group} of $H$ if $G$ is perfect, meaning the derived group $G'$ of $G$ equals $G$. Schur~\cite{Schur} demonstrated that a simple group $T$ possesses a universal covering group $G$, that is, every covering group of $T$ is a homomorphic image of $G$, and the center $Z(G)$ is called the {\em Schur multiplier} of $T$, denoted by $\Mult(T)$. The following result was presented in~\cite[Lemma 2.11]{PLHL}.

\begin{proposition}\label{mult}
Let $N$ be a group of order prime or prime-square, and let $T$ be a nonabelian simple group. Then every extension $G=N.T$ is a central extension. Furthermore, $G=NG'$ and $G'= M.T$ with $M\leq N\cap\Mult(T)$.
\end{proposition}

Let a group $G$ act on two sets $\Omega$ and $\Delta$. The two actions are termed {\em equivalent} if there exists a bijection $\lambda:\Omega\rightarrow\Delta$ such that $(\alpha^x)^\lambda=(\alpha^\lambda)^x$ for all $\alpha\in\Omega$ and $x\in G$. It is known that two transitive actions of $G$ on $\Omega$ and $\Delta$ are equivalent if and only if the stabilizers $G_\alpha$ and $G_\delta$ of $\alpha\in\Omega$ and $\delta\in\Delta$ in $G$ are conjugate in $G$ (see~\cite[Lemma 6.1B]{Dixon}).

Let $G$ act transitively on $\Omega$. A non-empty subset $\Delta$ of $\Omega$ is termed a {\em block} if for any $g\in G$, either $\Delta^g\cap \Delta=\emptyset$ or $\Delta^g=\Delta$, and a partition $\{\Omega_1,\cdots,\Omega_k\}$ of $\Omega$ is defined as {\em $G$-invariant} if every element of $G$ maps each $\Omega_i$ to some $\Omega_j$, where $1\leq i,j\leq k$.
Evidently, $\Omega$ and the singletons $\{\alpha\}$ ($\alpha\in\Omega$) are blocks of $G$, with such a block denoted as {\em trivial}, while any other block is referred to as {\em non-trivial}.
A transitive permutation group $G$ on $\Omega$ is called {\em primitive} if it has only trivial blocks in $\Omega$, and {\em quasiprimitive} if every non-trivial normal subgroup of $G$ is transitive on $\Omega$.
It is evident that a primitive group is quasiprimitive, but the converse does not hold true.

According to \cite[Corollary 3.4]{Li03}, a quasiprimitive permutation group containing a regular cyclic subgroup is primitive, and the classification of such permutation groups was independently achieved in~\cite{Jon} and~\cite{Li03}. In Tables~\ref{c-group} and~\ref{d-group} below,
the columns labeled ``$3$-tran" or ``$\#$Act" provide information on the $3$-transitivity and the number of inequivalent actions respectively. These details can be derived from the established classification of $2$-transitive permutation groups (refer to~\cite{LSS}),  \cite[Tables 1 and 2]{JT}, or~\cite[Theorems 3.2 and 3.3]{Pan}).

\begin{proposition}\label{Qc-group}
Let $G$ be a quasiprimitive permutation group on a set $\Ome$ such that $G$ contains a regular cyclic subgroup $H$ of degree $n$.
Let $\a\in\Ome$. Then either $n=p$ and $G\leq\AGL(1,p)$, or $G$ is $2$-transitive and $(G,G_\a,n)$ is listed in Table~$\ref{c-group}$.
\end{proposition}

\begin{table}[!htb]
\caption{Quasiprimitive permutation groups with a regular cyclic subgroup}
\label{c-group}
\begin{center}
\begin{tabular}{l|l|l|l|l|l} \hline
  $G$ & $G_\a$ & $n$ & Conditions & $3$-tran & $\#$Act \\ \hline
  $\A_n$ & $\A_{n-1}$ & $n$ & $n\geq 5$ is odd & Yes & 2 iff $n=6$ \\
  $\Sy_n$ & $\Sy_{n-1}$ & $n$ & $n\geq 4$ & Yes & 2 iff $n=6$ \\
  $\PGL(2,q).o$ & $[q]{:}\GL(1,q).o$ & $q+1$ & $o\leq\PGammaL(2,q)/\PGL(2,q)$ & Yes & $1$ \\
  $\PGL(d,q).o$ & $[q^{d-1}]{:}\GL(d-1,q).o$ & $\frac{q^d-1}{q-1}$ & $o\leq\PGammaL(d,q)/\PGL(d,q)$ & No & 2 \\
  $(d\geq3)$   &  &  &  &  \\
  $\PGammaL(2,8)$ & $\ZZ^3_2:\ZZ_7:\ZZ_3$ & $9$ &  & Yes &  $1$ \\
  $\PSL(2,11)$ & $\A_5$ & $11$ &  & No & 2\\
  $\M_{11}$ & $\M_{10}$ & $11$ &  & Yes & $1$  \\
  $\M_{23}$ & $\M_{22}$ & $22$ &  & Yes & $1$ \\ \hline
\end{tabular}
\end{center}
\end{table}

Let $V(d,q)$ be a linear space of dimension $d$ over the field $\F_q$ of order a prime-power $q$, and let $\PP_i$ be the set of all $i$-dimension subspaces of $V(d,q)$ for every $1\leq i\leq n$.
By Proposition~\ref{Qc-group}, $\PGL(d,q).o$ ($d\geq 3$) has two non-equivalent actions, which can viewed as its natural actions on $\PP_1$ and  $\PP_{n-1}$, respectively. Furthermore, the transpose inverse map of $\PSL(d,q)$ is an outer automorphism interchanging $\PP_1$ and $\PP_{n-1}$ (see~\cite[P .49]{Wilson}).

The quasiprimitive permutation groups containing a regular dihedral subgroup were classified in~\cite[Theorem 3.3]{SLZ}.

\begin{proposition}\label{Qd-group}
Let $G$ be a quasiprimitive permutation group on a set $\Ome$ such that $G$ contains a regular dihedral subgroup $H$.
Let $\a\in\Ome$. Then $G$ is $2$-transitive on $\Ome$ and $(G,G_\a,H)$ is listed in Table~$\ref{d-group}$.
\end{proposition}

\begin{table}[!htb]
\caption{Quasiprimitive permutation groups with a regular dihedral subgroup}
\label{d-group}
\begin{center}
\begin{tabular}{l|l|l|l|l|l} \hline
  $G$ & $G_\a$ & $H$ & Conditions & $3$-tran & $\#$Act\\ \hline
  $\A_4$ & $\ZZ_3$ & $\D_4$ &  & No  &  $1$ \\
  $\Sy_4$ & $\Sy_3$ & $\D_4$ &  & Yes  &  $1$ \\
  $\AGL(3,2)$ & $\GL(3,2)$ & $\D_8$ &  & Yes &  $1$ \\
  $\AGL(4,2)$ & $\GL(4,2)$ & $\D_{16}$ &  & Yes & $1$  \\
  $\ZZ_2^4:\A_7$ & $\A_7$ & $\D_{16}$ &  & Yes &  $1$ \\
  $\ZZ_2^4:\Sy_6$ & $\Sy_6$ & $\D_{16}$ &  & No &  $1$ \\
  $\ZZ_2^4:\A_6$ & $\A_6$ & $\D_{16}$ &  & No &  $1$ \\
  $\ZZ_2^4:\Sy_5$ & $\Sy_5$ & $\D_{16}$ &  & No &  $1$ \\
  $\ZZ_2^4:\GammaL(2,4)$ & $\GammaL(2,4)$ & $\D_{16}$ &  & No &  $1$ \\
  $\M_{12}$ & $\M_{11}$ & $\D_{12}$ &  & Yes & 2  \\
  $\M_{22}.\ZZ_2$ & $\PSL(3,4).\ZZ_2$ & $\D_{22}$ &  & Yes & $1$ \\
  $\M_{24}$ & $\M_{23}$ & $\D_{24}$ &  & Yes & $1$ \\
  $\Sy_{2n}$ & $\Sy_{2n-1}$ & $\D_{2n}$ &  & Yes  & $2$ iff $n=3$  \\
  $\A_{4n}$ & $\A_{4n-1}$ & $\D_{4n}$ &  & Yes  & $1$  \\
  $\PSL(2,r^f).o$ & $\ZZ_r^f:\ZZ_{\frac{r^f-1}{2}}.o$ & $\D_{r^f+1}$ & $r^f\equiv3\pmod{4}\text{~and~}$ & $3$-transitive & $1$ \\
   &  &  & $o\leq\ZZ_2\times\ZZ_f$ &  iff $\ZZ_2\leq o$  & \\
  $\PGL(2,r^f).\ZZ_e$ & $\ZZ_r^f:\ZZ_{r^f-1}.\ZZ_e$ & $\D_{r^f+1}$ & $r^f\equiv1\pmod{4},e\mid f$ & Yes &  $1$ \\
  \hline
\end{tabular}
\end{center}
\end{table}

A transitive permutation group is called {\em biquasiprimitive} if each of its nontrivial normal subgroups has at most two orbits and at least one has exactly two orbits. A transitive permutation group $G$ on $\Omega$ is called {\em biprimitive} if $G$ has an invariant partition $\{\Omega_1, \Omega_2\}$ such that $G_{\Omega_1}=G_{\Omega_2}$ is primitive on $\Omega_1$ and $\Omega_2$. Clearly, biprimitive group is biquasiprimitive, but the converse is not true. By Pan et al.~\cite[Theorem 1.1]{PYZH}, we have the following result.

\begin{proposition}\label{biquasi-d}
Each biquasiprimitive permutation group containing a regular dihedral subgroup is biprimitive.
\end{proposition}

It is well-known that $\GL(d,q)$ has a cyclic group of order $q^d-1$, the so-called {\em Singer cycle} of $\GL(d,q)$, which induces an cyclic group of order $(q^d-1)/(q-1)$ in $\PGL(d,q)$, called a {\em Singer cycle} of $\PGL(d,q)$.

\begin{lemma}\label{proper-PSL}
\begin{enumerate}[\rm (1)]
  \item Let $\la x\ra$ be a Singer cycle of $\PGL(d,q)$ with $d\geq 3$, and let $\la x\ra< G<\PGammaL(d,q)$ with $|G:\la x\ra|=2$. Then $G$ cannot be cyclic or dihedral.
  \item Let $\PSL(2,q).o\leq \PGammaL(2,q)$ with $q=r^f\equiv3\pmod{4}$ for a prime $r$, and $o\leq\ZZ_2\times\ZZ_f$ does not contain the diagonal automorphism of $\PSL(2,q)$. Then every dihedral group of order $q+1$ in $\PSL(2,q).o$ is a subgroup of $\PSL(2,q)$.
\end{enumerate}
\end{lemma}

\proof Note that $\langle x \rangle \cong \mathbb{Z}_{(q^d-1)/(q-1)}$. Since $|G:\langle x \rangle|=2$, we have $\langle x \rangle \trianglelefteq G$, and so $G \leq N_{\PGammaL(d,q)}(\langle x \rangle)$.
By~\cite[P. 187, 7.3 Satz]{Hupp}, $N_{\PGammaL(d,q)}(\langle x \rangle)=\langle x \rangle : \langle \delta \rangle$ with $\langle \delta \rangle \cong \mathbb{Z}_d$ and $x^\delta=x^q$.
Assume that $G$ is cyclic or dihedral.
Then there is an integer $0 \leq i < d$ such that $x^{\delta^i}=x^{q^i}=x^{\pm 1}$, that is, $x^{q^{i\pm 1}}=1$.
This implies that the order of $x$ is a divisor of $q^i\pm 1$, that is, $(q^d-1)/(q-1)=q^{d-1}+\cdots+q+1$ divides $q^i\pm 1$, which leads to $d=2$ and $i=1$, contradicting the assumption $d \geq 3$. This completes the proof of part (1).

To prove part (2), let $L\leq \PSL(2,q).o$ with $L\cong \D_{q+1}$.
Suppose that $L\not\leq  \PSL(2,q)$.
Since $o\leq\ZZ_2\times\ZZ_f$ does not contain the diagonal automorphism of $\PSL(2,q)$, we have that $o$ is cyclic. Note that $L\cap \PSL(2,q)\unlhd L$ and $L/(L\cap \PSL(2,q))\cong \PSL(2,q) L/\PSL(2,q)\leq o$. This implies $L\cap \PSL(2,q)\not=1$, and since $L$ is dihedral and $q+1\geq 8$ ($q\equiv3\pmod{4}$), we have $|L:L\cap \PSL(2,q)|=2$, forcing $2\mid |o|$.
Thus, $f$ is even as $|o|\mid f$. It follows that $q\equiv 1\pmod{4}$ from $q=r^f$, contradicting to $q\equiv3\pmod{4}$. This yields that $L\leq  \PSL(2,q)$, as part (2).
\qed

\subsection{Graph theoretical results}
Let $q=p^f$ be a prime power such that $q\equiv1\pmod{4}$, and let $\F_q$ be the finite field of order $q$. The {\em Paley graph} $\P(q)$ is defined to be the graph with vertex set $\F_q$, and two vertices $u,v$ are adjacent if and only if $u-v$ is a nonzero square in $\F_q$.
This graph was first defined by Paley~\cite{Paley}, and $\P(q)=\Cay(\F_q^+,S)$ with $S=\la\lambda^2\ra$, where $\F_q^+$ is the additive group of $\F_q$ and $\lambda$ be a primitive element of $\F_q$. Furthermore,  $\P(q)$ has valency $(q-1)/2$ and diameter $2$, and $\P(q)$ has girth $3$ for $q>5$ and $\P(5)\cong\C_5$.

A Cayley graph $\Cay(G,S)$ is called {\em circulant} if $G$ is a cyclic group. The 2-arc transitive and 2-distance transitive circulants were classified in~\cite{ACX,CJL}, respectively. We summarize those results as follows:

\begin{proposition}\label{cir}
Let $\Ga$ be a connected circulant. Then the following holds:
\begin{enumerate}[\rm (1)]
  \item If $\Ga$ is $2$-arc transitive, then $\Ga$ is one of the following graphs: $\K_n$ with $n\geq 1$, $\C_n$ with $n\geq 4$, $\K_{\frac{n}{2},\frac{n}{2}}$ with $n\geq 6$, $\K_{\frac{n}{2},\frac{n}{2}}-\frac{n}{2}\K_2$ with $\frac{n}{2}\geq5$ odd.
  \item If $\Ga$ is $2$-distance transitive but not $2$-arc transitive, then $\Ga$ is isomorphism to either $\K_{m[b]}$ for some $m\geq3$ and $b\geq2$, or the Paley graph $\P(p)$ with $p$ a prime and $p\equiv1\pmod{4}$.
\end{enumerate}
\end{proposition}

The next result is quoted from~\cite[Lemma 5]{JWZ}.

\begin{proposition}\label{girth4}
Let $\Ga$ be a $2$-distance transitive graph of girth $4$ and valency $r\geq 3$.
If $|\Ga_2(u)|=r$ for some $u\in V(\Ga)$, then $\Ga\cong\K_{r+1,r+1}-(r+1)\K_2$.
\end{proposition}

Let $X$ be a group of automorphisms of a graph $\Ga$. For a normal subgroup $N$ of $X$, the {\em quotient graph} $\Ga_N$ of $\Ga$ with respect to $N$ is defined to be the graph with the orbit set of $N$ as vertex set and with two distinct orbits $B,C$ of $N$ adjacent in $\Ga_N$ if there is an edge of $\Ga$ between $B$ and $C$.
The graph $\Ga$ is called a {\em cover} or an {\em $N$-cover} of $\Ga_N$ if every vertex $u\in V(\Ga)$ has the same valency as $u^N$ in $\Ga_N$, where $u^N$ is the orbit of $N$ containing $u$ in $V(\Ga)$.

A graph $\Ga$ is called {\em locally $(G,s)$-distance transitive} with $s\geq 1$, if for every vertex $u\in V(\Ga)$, $G_u$ is transitive on $\Ga_i(u)$ for all $i\leq s$.
The next result illustrates that a connected locally $(G,2)$-distance transitive graph cannot be a cover of a complete multipartite graph.

\begin{proposition}[{\cite[Proposition 4.2]{DGLP12}}]\label{Kmb}
Let $\Ga$ be a connected locally $(G,s)$-distance transitive graph with $s\geq 2$. Then there exists no nontrivial $N\unlhd G$ such that $\Ga$ is a cover of $\Ga_N$ with $\Ga_N\cong\K_{m[b]}$, where $m\geq3$ and $b\geq 2$.
\end{proposition}

The following proposition gives a reduction for studying locally $(G,s)$-distance transitive graphs.

\begin{proposition}[{\cite[Lemma 5.3]{DGLP12}}]  \label{redu}
Let $\Ga$ be a connected locally $(G,s)$-distance transitive graph with $s\geq2$. Let $1\neq N\unlhd G$ be intransitive on $V(\Ga)$, and let $\BB$ be the set of $N$-orbits on $V(\Ga)$. Then one of the following holds:
\begin{enumerate}[\rm (1)]
  \item $|\BB|=2$;
  \item $\Ga$ is bipartite, $\Ga_N\cong\K_{1,r}$ with $r\geq2$ and $G$ is intransitive on $V(\Ga)$;
  \item $s=2$, $\Ga\cong\K_{m[b]}$ and $\Ga_N\cong\K_m$, where $m\geq3$ and $b\geq 2$;
  \item $N$ is semiregular on $V(\Ga)$, $\Ga$ is a cover of $\Ga_N$, $|V(\Ga_N)|<|V(\Ga)|$ and $\Ga_N$ is locally $(G/N,s')$-distance transitive, where $s'=\min\{s,\diam(\Ga_N)\}$.
\end{enumerate}
\end{proposition}

The following are well-known results, of which the first two parts follow from
\cite[Lemma 2.5]{LP} or \cite[Theorem 4.1]{Praeger93}, and the third part follows from Proposition~\ref{redu}.

\begin{proposition}\label{kernel}
Let $\Ga$ be a connected graph and let $X\leq\Aut(\Ga)$. Assume that $N\unlhd X$ such that $\Ga$ is an $N$-cover of the quotient graph $\Ga_N$. Then the following holds.
\begin{enumerate}[\rm (1)]
  \item $N$ is the kernel of $X$ acting on $V(\Ga_N)$, and $N$ is semiregular on $V(\Ga)$.
  \item $\Ga$ is $X$-arc transitive or $(X,2)$-arc transitive if and only if $\Ga_N$ is $X/N$-arc transitive or $(X/N,2)$-arc transitive, respectively.
  \item If $\Ga$ is $(X,2)$-distance transitive then $\Ga_N$ is $(X/N,2)$-distance transitive.
\end{enumerate}
\end{proposition}

By Corr et al.~\cite[Lemma 3.3]{CJS}, a subgroup of $\Aut(\K_{n,n})$ is $2$-distance transitive if and only if it is $2$-arc transitive on $\K_{n,n}$.
The proof of this result relies on the classification of $2$-transitive permutation groups and so on the classification of finite simple groups.
In the subsequent lemma, we provide a concise proof based on an elementary group result. The equivalence of $2$-distance transitivity and $2$-arc transitivity for a subgroup of automorphisms of a graph also holds for $\B(H_{11})$ or $\B'(H_{11})$, which are the incidence or non-incidence graphs of the Hadamard design on $11$ points, respectively.

\begin{lemma}\label{bipart}
Let $\Ga=\K_{n,n}$, $\B(H_{11})$ or $\B'(H_{11})$, and let $X\leq\Aut(\Ga)$. Then $\Ga$ is $(X,2)$-distance transitive if and only if $\Ga$ is $(X,2)$-arc transitive.
\end{lemma}

\proof To complete the proof, it is sufficient to demonstrate the necessity, which entails proving that $\Ga$ is $(X,2)$-arc transitive given that $\Ga$ is $(X,2)$-distance transitive. Let $(u,v,w)$ be a $2$-arc of $\Ga$. As $\Ga$ is bipartite, it follows that $v\in \Ga(u)$ and $w\in \Ga_2(u)$. Moreover, since $\Ga$ is $(X,2)$-distance transitive, $X_u$ exhibits transitivity on both $\Ga(u)$ and $\Ga_2(u)$.

Assume $\Ga=\K_{n,n}$. Then $|\Ga(u)|=n$ and $\Ga_2(u)=n-1$. The lemma trivially holds for $n=1$, so we will proceed with the assumption that $n\geq 2$.
By utilizing the orbit-stabilizer theorem~\cite[Theorem 1.4A]{Dixon}, we find that $|X_u:X_{uv}|=n$ and $|X_u:X_{uw}|=n-1$. Since $(n,n-1)=1$, based on \cite[P. 8, 2.13 Hilfssatz]{Hupp}, we have $|X_u:X_{uvw}|=|X_u:X_{uv}|\cdot|X_u:X_{uw}|$, indicating $|X_u|=|X_{uv}||X_{uw}|/|X_{uvw}|$.
Furthermore, as $|X_{uv}X_{uw}|=|X_{uv}||X_{uw}|/|X_{uv}\cap X_{uw}|=|X_{uv}||X_{uw}|/|X_{uvw}|$, we obtain $X_u=X_{uw}X_{uv}$. Thus, $X_{uv}$ can map $w\in \Ga_2(u)$ to any vertex in $\Ga_2(u)$ because $X_u$ is transitive on $\Ga_2(u)$. Since $X$ is arc transitive, $X$ is $2$-arc transitive.

Assume $\Ga=\B(H_{11})$ or $\B'(H_{11})$.
In these cases, we have $\Aut(\Gamma)\cong\PGL(2,11)$, and $\Gamma$ is a $2$-arc transitive graph of order $22$ as stated in~\cite[Example 1.1]{Pan}.
Notice that all subgroups of $\PSL(2,11)$ and $\PGL(2,11)$ are listed in~\cite[P. 7]{Atlas}.
Since $\Gamma$ is $X$-vertex transitive, $X$ must have a subgroup of index $22$, implying  $X\not\cong\PSL(2,11)$ because $\PSL(2,11)$ has no subgroup of index $22$.
First let $\Ga=\B'(H_{11})$. Then $\Ga$ has girth $4$ and valency $6$, and hence $6\times 22$ divides $|X|$. It follows that $X\cong\PGL(2,11)$, and $\Ga$ is $(X,2)$-arc transitive, as required.
Now let $\Ga=\B(H_{11})$. Then $\Ga$ has valency $5$ and $110$ divides $|X|$, implying  $X\cong\ZZ_{11}:\ZZ_{10}$ or $\PGL(2,11)$.
If $X\cong\ZZ_{11}:\ZZ_{10}$, then $X_u\cong\ZZ_5$, and so $|\Ga(u)|=|\Ga_2(u)|=5$.
By Proposition~\ref{girth4}, $\Ga\cong\K_{6,6}-6\K_2$, a contradiction.
Thus $X=\PGL(2,11)$ and $\Ga$ is $(X,2)$-arc transitive, thereby concluding the proof.
\qed

For a graph $\Gamma$ and subsets $U,W\subseteq V(\Gamma)$, denote by $[U]_{\Gamma}$ the induced subgraph of $U$ in $\Gamma$, and by $[U,W]_{\Gamma}$ the subgraph of $\Gamma$ with $V([U,W]_{\Gamma})=U\cup W$ and $E([U,W]_{\Gamma})=\{\{u,w\}\ |\ u\in U,w\in W\mbox{ and } \{u,w\}\in E(\Gamma)\}$.
Clearly, if $U\cap W=\emptyset$ then $[U,W]_\Gamma$ is bipartite. If all vertices of $[U,W]_\Gamma$ have same valency, we denote it by $v([U,W]_\Gamma)$ the same valency, and similarly, if all vertices of $[U]_{\Gamma}$ have same valency, we denote it by $v([U]_{\Gamma})$. Furthermore, we will omit the subscript $\Gamma$ in $[U]_\Gamma$ and $[U,W]_\Gamma$ if there is no confusion.
To end the section, we consider a special family of connected $2$-distance transitive bipartite graphs.

\begin{lemma}\label{bipartitegraphs}
Let $\Ga$ be a connected $(X,2)$-distance transitive bipartite graph of valency $r\geq 3$, where $X\leq\Aut(\Ga)$. For $u\in V(\Gamma)$, assume that $X_u$ is regular on the neighbourhood $\Gamma(u)$ of $u$ in $\Gamma$. Then $\Gamma\cong \K_{r,r}$, or $\K_{r+1,r+1}-(r+1)\K_2$.
\end{lemma}

\proof Let $\Delta_1$ and $\Delta_2$ be the bipartite sets of $\Ga$ with $u\in \Delta_1$. Then $\Gamma(u)\subseteq \Delta_2$ and $\Gamma_2(u)\subseteq \Delta_1$. Since $X_u$ is regular on $\Gamma(u)$, we have $|X_u|=|\Gamma(u)|=r$.

Consider the induced subgroup of $[\Gamma(u)\cup \Gamma_2(u)]$ of $\Gamma(u)\cup \Gamma_2(u)$ in $\Gamma$. Clearly, $[\Gamma(u)\cup \Gamma_2(u)]$ is a bipartite graph with partite sets $\Gamma(u)$ and $\Gamma_2(u)$. Every vertex in $\Gamma(u)$ has all neighbours in $\Gamma_2(u)$ except $u$, thereby having $r-1$ neighbors. This implies that the number of edges in $[\Gamma(u)\cup \Gamma_2(u)]$  is $r(r-1)=|\Gamma(u)|(r-1)$. On the other hand, since $X_u$ is transitive on $\Gamma_2(u)$, every vertex in $\Gamma_2(u)$ has the same number of neighbours, say $\ell$, in $\Gamma(u)$,  leading to the number of edges in $[\Gamma(u)\cup \Gamma_2(u)]$ being $\ell |\Gamma_2(u)|$. It follows that $\ell |\Gamma_2(u)|=r(r-1)$. Note that $\ell\leq r$. Since $|X_u|=r$ and $X_u$ is transitive on $\Gamma_2(u)$, we have $|\Gamma_2(u)|\leq r$. Then $\ell |\Gamma_2(u)|=r(r-1)$ implies that either $\ell=r$ and $|\Gamma_2(u)|=r-1$, or $|\Gamma_2(u)|=r$ and $\ell=r-1$.
For the former, the connectivity of $\Gamma$ suggests $\Gamma=\K_{r,r}$.
For the latter, $[\Gamma(u)\cup \Gamma_2(u)]=\K_{r,r}-r\K_2$ with $|\Gamma(u)|=|\Gamma_2(u)|=r$, and since $r\geq 3$, $\Gamma$ has girth $4$ (a $4$-cycle consists of $u$, two vertices in $\Gamma(u)$ and one vertex in $\Gamma_2(u)$). By Proposition~\ref{girth4}, we conclude that $\Gamma=\K_{r+1,r+1}-(r+1)\K_2$. \qed

\section{$2$-distance transitive covers of graphs}\label{properties}

In this section, we delve into the existence of $2$-distance transitive covers of given graphs under specific conditions.

To begin, we aim to demonstrate that a connected graph cannot serve as a $2$-distance transitive cover of the Paley graph $\P(p)$ for a prime $p>5$.

\begin{lemma}\label{Paley}
Let $\Ga$ be a connected $(X,2)$-distance transitive graph with $X\leq\Aut(\Ga)$.
Assume that $1\neq N\unlhd X$ is intransitive on $V(\Ga)$.
Then $\Ga$ cannot be an $N$-cover of $\P(p)$, where $p$ is a prime with $p\equiv1\pmod{4}$ and $p> 5$.
\end{lemma}

\proof Suppose to the contrary that $\Ga$ is an $N$-cover of $\Sig:=\Ga_N\cong\P(p)$. By Proposition~\ref{kernel}, $N$ is the kernel of $X$ acting on the orbit set $V(\Sig)$ of $N$ in $V(\Ga)$ and $X/N\leq \Aut(\Sig)$. Furthermore, $\Sig$ is $(X/N,2)$-distance transitive. It is noted that $p\geq 13$ as $p\equiv1\pmod{4}$ and $p> 5$.
Let $B_0\in V(\Sigma)$ with $b_0\in B_0$. Since $\Sigma\cong\P(p)$ has valency $(p-1)/2$ and diameter $2$, we may let
$$\Sigma(B_0)=\{B_1,B_2,\ldots,B_{(p-1)/2}\} \text{ and } \Sigma_2(B_0)=\{C_1,C_2,\ldots,C_{(p-1)/2}\}.$$
Therefore, $V(\Gamma_N)=\{B_0\}\cup \Sigma(B_0)\cup \Sigma_2(B_0)$.
By~\cite[P. 221]{GR}, $|\Sigma(B_1)\cap\Sigma(B_0)|=(p-5)/4$, hence $|\Sigma(B_1)\cap\Sigma_2(B_0)|=(p-1)/4$. This implies that
$v([\Sigma(B_0)]_\Sigma)=(p-5)/4$, and $v([\Sigma(B_0),\Sigma_2(B_0)]_\Sigma)=(p-1)/4$.

The orbit set of $N$ on $V(\Gamma)$ is $\{B_0,B_1,\cdots,B_{(p-1)/2},C_1,C_2,\cdots,C_{(p-1)/2}\}$, and $X$ has a natural action on the orbit set. Since $\Gamma$ is a cover of $\Sig$ and $b_0\in B_0$, $\Gamma$ has valency $(p-1)/2$ and we may write
$$\Gamma(b_0)=\{b_1,b_2,\cdots,b_{(p-1)/2}\} \mbox{ with } b_i\in B_i \mbox{ for }1\leq i\leq (p-1)/2.$$
Assume that $(b_0,b_1,c_1)$ is a $2$-arc in $\Gamma$ for some $c_1\in C_1$. Then $d_\Gamma (b_0,c_1)=2$, where $d_\Gamma (b_0,c_1)$ is the distance between $b_0$ and $c_1$ in $\Gamma$. Since $\Gamma$ is $(X,2)$-distance transitive, $X_{b_0}$ is transitive on both $\Gamma(b_0)$ and $\Gamma_2(b_0)$.

Take an edge $\{b_i,b_j\}$ in the induced subgraph $[\Gamma(b_0)]_\Ga$. Then $\{B_i,B_j\}$ is an edge in $[\Sig(B_0)]_\Sig$ as $\Gamma$ is an $N$-cover of $\Sig$. On the other hand, take an edge $\{B_i,B_j\}$ in $[\Sig(B_0)]_\Sig$. Then $b_i$ is adjacent to some vertices in $B_j$, say $b$. If $b\not=b_j$ then $d_\Gamma (b_0,b)=2$. Since $\Gamma$ is $(X,2)$-distance transitive, there is $\alpha\in X_{b_0}$ such that $b^\alpha=c_1$ as $d_\Gamma (b_0,c_1)=2$, and hence $\alpha N\in X/N\leq \Aut(\Sig)$ fixes $B_0$ and maps $B_j$ to $C_1$, which is impossible because $B_0$ is adjacent to $B_j$ in $\Sig$ but not to $C_1$. This implies that $\{b_i,b_j\}$ is an edge in $[\Gamma(b_0)]_\Gamma$. It follows that $[\Gamma(b_0)]_\Ga\cong [\Sig(B_0)]_\Sig$, and hence $v([\Gamma(b_0)]_\Ga)=v([\Sig(B_0)]_\Sig)=(p-5)/4$. In particular, every vertex in $\Gamma(b_0)$ has valency  $(p-1)/4$ in the bipartite subgraph $[\Gamma(b_0),\Gamma_2(b_0)]_\Ga$ because $\Ga$ has valency $(p-1)/2$.

Since $\Gamma$ is a cover of $\Sig$, $|\Gamma_2(b_0)|\geq |\Sig_2(B_0)|=(p-1)/2$ and $X_{b_0}\cong (X/N)_{B_0}$. By \cite[Theorem 7.1(1)]{Pei}, $\Aut(\P(p))\cong \ZZ_p:\ZZ_{(p-1)/2}$ and hence $|X_{b_0}|=|(X/N)_{B_0}|\leq (p-1)/2$. Since $X_{b_0}$ is transitive on $\Gamma_2(b_0)$, we have $|\Gamma_2(b_0)|\leq (p-1)/2$. It follows that $|\Gamma_2(b_0)|=(p-1)/2=|\Gamma(b_0)|$. Since every vertex in $\Gamma(b_0)$ has valency  $(p-1)/4$ in $[\Gamma(b_0),\Gamma_2(b_0)]_\Ga$, we have  $v([\Gamma(b_0),\Gamma_2(b_0)]_\Ga)=(p-1)/4$. Write $$\Gamma_2(b_0)=\{c_1,c_2,\cdots,c_{(p-1)/2}\} \mbox{ with }c_i\in C_i \mbox{ for } 1\leq i\leq (p-1)/2.$$
Then $|\Gamma(b_0)\cap \Gamma(c_1)|=(p-1)/4$.
Note that $|\Gamma(b_0)|+|\Gamma_2(b_0)|=p-1$. Since $N\not=1$, we have $|V(\Gamma)|>p$, and hence $\Gamma_3(b_0)\not=\emptyset$. Since $X_{b_0}$ acts transitively on $\Gamma_2(b_0)$, we have $\Gamma(c_1)\cap \Gamma_3(b_0)\not=\emptyset$.
Let $f\in \Gamma(c_1)\cap \Gamma_3(b_0)$. Then $(b_0,b_1,c_1,f)$ is a $3$-arc and $d_\Gamma (b_0,f)=3$.
It follows that $d_\Gamma (b_1,f)=2$, and the  $(X,2)$-distance transitivity of $\Gamma$ gives rise to $|\Gamma(b_1)\cap \Gamma(f)|=|\Gamma(b_0)\cap \Gamma(c_1)|=(p-1)/4$, forcing $|\Gamma_2(b_0)\cap \Gamma(f)|\geq (p-1)/4$.
Since $f\in \Gamma(c_1)$ and $c_1\in C_1$, we have $f\in B_i$ or $C_i$ for some $1\leq i\leq (p-1)/2$.

Note that $X_{b_0}$ acts transitively on both $\Sig(B_0)$ and $\Sig_2(B_0)$ as it is transitive on $\Gamma(b_0)$ and $\Gamma_2(b_0)$. Then for every  $1\leq i\leq (p-1)/2$, $B_i$ or $C_i$ contains an element, say $f_i$, such that $d_\Gamma (b_0,f_i)=3$ and $|\Gamma_2(b_0)\cap \Gamma(f_i)|\geq (p-1)/4$, because $d_\Gamma (b_0,f)=3$ and $|\Gamma_2(b_0)\cap \Gamma(f)|\geq (p-1)/4$. Set $$F=\{f_1,f_2,\cdots,f_{(p-1)/2}\}.$$
Then $F\subseteq \Gamma_3(b_0)$. Since $|\Gamma_2(b_0)\cap \Gamma(f_i)|\geq (p-1)/4$,
every vertex in $F$ has valency at least $(p-1)/4$ in the bipartite subgraph $[\Gamma_2(b_0),F]_\Ga$,
and hence $|E([\Gamma_2(b_0),F]_\Ga)|\geq (p-1)/2\cdot (p-1)/4$.
Since $v([\Gamma(b_0),\Gamma_2(b_0)]_\Ga)=(p-1)/4$,
every vertex in $\Gamma_2(b_0)$ has valency at most $(p-1)/4$ in the bipartite subgraph $[\Gamma_2(b_0),F]_\Ga$, and hence $|E([\Gamma_2(b_0),F]_\Ga)|\leq (p-1)/2\cdot (p-1)/4$.
It follows that $|E([\Gamma_2(b_0),F]_\Ga)|=(p-1)/2\cdot (p-1)/4$, forcing  $v([\Gamma_2(b_0),F]_\Gamma)=(p-1)/4$. This also implies that $F=\Gamma_3(b_0)$ and $v([\Gamma_2(b_0)])=0$.

For an edge $e$ in $\Gamma$, denote by $tri(e)$ the number of triangles passing through $e$ in $\Gamma$. Since $v([\Gamma(b_0)]_\Ga)=(p-5)/4$, we have $tri(\{b_0,b_1\})=(p-5)/4$, and the arc transitivity of $\Ga$ implies that  $tri(e)=(p-5)/4$. Take $\{b_1,c_i\}\in E([\Gamma(b_0),\Gamma_2(b_0)]_\Ga)$.  Then $tri(\{b_1,c_i\})=(p-5)/4$, and since $v([\Gamma_2(b_0)])=0$, $c_i$ is adjacent to every neighbour of $b_1$ in $[\Gamma(b_0)]_\Ga$ as $v([\Gamma(b_0)]_\Ga)=(p-5)/4$. Then   $[\Gamma(b_0),\Gamma_2(b_0)]_\Ga\cong 2\K_{(p-1)/4,(p-1)/4}$ by  the arbitrariness of $c_i$. In particular, for every edge $\{b_1,b_i\}\in E([\Gamma(b_0)]_\Ga)$, $b_1$ and $b_i$ belong to the same component of $[\Gamma(b_0),\Gamma_2(b_0)]_\Ga$ isomorphic to $\K_{(p-1)/4,(p-1)/4}$. It follows $tri(\{b_1,b_i\})=(p-1)/4+1$, contradicting to $tri(e)=(p-5)/4$. This completes the proof.
\qed

The forthcoming result illustrates that a connected graph cannot act as a $2$-distance transitive cover of a particular bipartite graph under specific conditions.

\begin{lemma}\label{2pcover}
 Let $\Ga$ be a connected $(X,2)$-distance transitive graph and let $1\not=N\unlhd X$ with $X/N$ solvable. Then $\Ga$ cannot be an $N$-cover of $\K_{p,p}-p\K_2$ for a prime $p\geq 5$.
\end{lemma}

\proof Let $\Sig=\K_{p,p}-p\K_2$ for a prime $p\geq 5$.
Assume, for the sake of contradiction, that $\Gamma$ is an $N$-cover of $\Sigma$. Then $\Gamma$ is a bipartite graph of valency $k=p-1$. By Proposition~\ref{kernel}, the kernel of $X$ on $V(\Gamma_N)$ is $N$, and hence $Y:=X/N\leq \Aut(\Gamma_N)$. We may let $\Ga_N=\Sigma$. Then $\Sigma$ is $(Y,2)$-distance transitive, and so $Y$ is arc transitive on $\Sigma$. Let $\Delta_1$ and $\Delta_2$ be the bipartite sets of $\Sigma$. Then $|\Delta_1|=|\Delta_2|=p$. Write $Y^*=Y_{\Delta_1}$, the subgroup of $Y$ fixing $\Delta_1$ setwise. Then $|Y:Y^*|=2$. Since $Y$ is solvable, $Y^*$ is solvable, and since $N\not=1$, we have $|V(\Gamma)|=|N||V(\Sigma)|=2p|N|>2p$.

Let $K$ denote the kernel of the action of $Y^*$ on $\Delta_1$. Then $K \unlhd Y^*$. If $K \neq 1$, then $K$ is transitive on $\Delta_2$ since $|\Delta_2|=p$. This would imply $\Sigma \cong \K_{p,p}$, contradicting the fact that $k=p-1$. Therefore, $Y^*$ is faithful on both $\Delta_1$ and $\Delta_2$.
According to \cite[Corollary 3.5B]{Dixon}, we have $Y^* \leq \text{AGL}(1,p)$, and we can express $Y^*=\ZZ_p:\ZZ_r$, where $r \mid (p-1)$.

For $u\in V(\Gamma)$, let $\overline{u}=u^N$ be the orbit of $u$ under $N$. Then $\overline{u}\in V(\Sigma)$. Since all subgroups of order $r$ of $Y^*=\ZZ_p: \ZZ_r$ are conjugate, there exist $\overline{u}\in\Delta_1$ and $\overline{v}\in\Delta_2$, where $u,v\in V(\Gamma)$, such that $\ZZ_r$ fixes $\overline{u}$ and $\overline{v}$ and semiregular on $\Delta_1\backslash\{\overline{u}\}$ and $\Delta_2\backslash\{\overline{v}\}$. Since $Y$ is arc transitive on $\Sigma$, $r$ is the valency of $\Sigma$ and so $r=p-1$. In particular, $Y_{\overline{u}}=Y^*_{\overline{u}}=\mz_{p-1}$. Since $\Gamma$ is a cover of $\Sigma$,  we have $X_u\cong\mz_{p-1}$, and since $\Gamma$ has valency $p-1$, $X_u$ is regular on the neighbourhood  $\Gamma(u)$ of $u$ in $\Gamma$. By Lemma~\ref{bipartitegraphs}, $\Gamma=\K_{p,p}-p\K_2$, contradicting to $|V(\Gamma)|>2p$. \qed

The subsequent result simplifies a $2$-distance cover of a graph into a combination of two $2$-distance covers.

\begin{lemma}\label{zp-cover}
Let $\Ga$ be a connected $(X,2)$-distance transitive graph, and an $N$-cover of $\Ga_N$ for $N\unlhd X$. Suppose $K\unlhd X$ and $K\leq N$. Then $\Ga$ is a $K$-cover of $\Ga_K$ and $\Ga_K$ is an $N/K$-cover of $\Ga_N$. Moreover, $\Ga_K$ is $(X/K,2)$-distance transitive and $\Ga_N$ is $(X/N,2)$-distance transitive.
\end{lemma}

\proof Consider the quotient graph $\Ga_K$ of $\Ga$ under $K$, and the quotient graph $(\Ga_K)_{N/K}$ of $\Ga_K$ under $N/K$. For every $u\in V(\Ga)$, if we identify the orbit $u^N$ of $u$ under $N$ with the orbit $(u^K)^{N/K}$ of $u^K\in V(\Ga_K)$ under $N/K$, it is easy to see that $\Ga_N=(\Ga_K)_{N/K}$. Since $\Ga$ is an $N$-cover of $\Ga_N$, it follows that $\Ga$, $\Ga_K$ and $\Ga_N$ have the same valency. Thus, $\Ga$ is a $K$-cover of $\Ga_K$ and $\Ga_K$ is an $N/K$-cover of $\Ga_N$. By Proposition~\ref{kernel}, $K$ is the kernel of $X$ acting on $V(\Ga_K)$ with $X/K\leq \Aut(\Ga_K)$, $\Ga_K$ is $(X/K,2)$-distance transitive and $\Ga_N$ is $(X/N,2)$-distance transitive.
\qed

\section{Proof of Theorem~\ref{Thm-1}}\label{proof1}

For convenience, throughout this section we set
$$\D_{2n}=\la x,y\mid x^n=y^2=1,x^y=x^{-1}\ra \text{~with~} n\geq2, \text{~and~} H=\la x\ra.$$

Initially, we characterize normal subgroups $N$ of $X$ for a $(X,2)$-distance transitive dihedrant that is an $N$-cover of the  quotient graph corresponding to $N$.

\begin{lemma}\label{base}
Let $\Ga=\Cay(\D_{2n},S)$ be a connected $(X,2)$-distance transitive graph with $n\geq2$ and $R(\D_{2n})\leq X\leq\Aut(\Ga)$.
Let $N$ be a normal subgroup of $X$ and let $\Ga$ be a $N$-cover of $\Ga_N$.
Then one of the following holds:
\begin{enumerate}[\rm (1)]
  \item $N<R(H)$ and $\Ga_N$ is a dihedrant on $R(\D_{2n})/N$;
  \item $|N:N\cap R(H)|=2$ and $\Ga_N$ is a circulant on $R(H)N/N$.  Moreover, there is a normal subgroup $M$ of $X$ such that $|N:M|=2$ or $4$.
\end{enumerate}
\end{lemma}

\proof Let $G=\D_{2n}$ with $n\geq2$. Let $V(\Gamma_N)=\{\Delta_1,\Delta_2,\cdots,\Delta_m\}$ be the set of orbits of $N$ on $V(\Ga)$. Then $|G|=|V(\Ga)|=2n\geq 4$. Since $\Ga$ is connected, it has valency at least $2$, and since $\Ga$ is a cover of $\Ga_N$, $\Ga_N$ has valency at least $2$, implying $m\geq 3$. By Proposition~\ref{kernel}, the kernel of $X$ on $\{\Delta_1,\Delta_2,\cdots,\Delta_m\}$ is $N$ that is semiregular, and $X/N\leq \Aut(\Gamma_N)$. Since $N\unlhd X$ and $R(G)\leq X$, $V(\Ga_N)$ is a complete imprimitive block system of $X$ on $V(\Ga)$, yielding $|\Delta_1|=\cdots=|\Delta_m|$, and in particular, $V(\Ga_N)$ is also a complete imprimitive block system of $R(G)$. Since $R(H)\unlhd R(G)$, $\{H,yH\}$ forms a complete imprimitive block system of $R(G)$ and $R(y)$ interchanges $H$ and $yH$.

First assume that $H$ or $yH$ contains a block in $V(\Ga_N)$. Since $R(y)$ interchanges $H$ and $yH$, both $H$ or $yH$ contain blocks in $V(\Ga_N)$, say $\Delta_1\subseteq H$ and $\Delta_2\subseteq yH$. Since $R(H)$ is transitive on $H$ and $yH$, both $H$ and $yH$ are union of some $\Delta_i$'s in $V(\Ga_N)$.

Note that $R(H)$ is regular on each of $H$ and $yH$. Since $\Delta_1$ is a block of $R(H)$, the block stabilizer $R(H)_{\Delta_1}$ of $\Delta_1$ in $R(H)$ is a subgroup of order $|\Delta_1|$ in $R(H)$, which acts regularly on $\Delta_1$. Since $R(H)$ is abelian, $R(H)_{\Delta_1}=R(h)^{-1}R(H)_{\Delta_1}R(h)=R(H)_{\Delta_1^h}$ for every $h\in H$, and hence $R(H)_{\Delta_1}$ fixes every $\Delta_i\subseteq H$ setwise. Similarly, the block stabilizer $R(H)_{\Delta_2}$ of $\Delta_2$ in $R(H)$ has order $|\Delta_2|$ and fixes every block $\Delta_j\subseteq yH$ setwise.
Since $R(H)$ is cyclic and $|\Delta_1|=|\Delta_2|$, we have $R(H)_{\Delta_1}=R(H)_{\Delta_2}$, and hence $R(H)_{\Delta_1}$ fixes every block in $V(\Ga_N)$. Since $N$ is the kernel of $X$ on $V(\Ga_N)$, we have $R(H)_{\Delta_1}\leq N$, and since $N$ is semiregular on $V(\Ga)$, we have $N_\a=1$ for $\a\in \Delta_1$. By the Frattini argument, $N=R(H)_{\Delta_1}N_\a=R(H)_{\Delta_1}$, and hence $N<R(H)$ as $m\geq 3$. Additionally, $\Ga_N$ is a dihedrant on $R(G)/N$.

Now we assume that neither $H$ nor $yH$ contains a block in $V(\Ga_N)$. Then for every $1\leq i\leq m$, $\Delta_i\cap H\not=\emptyset$ and $\Delta_i\cap yH\not=\emptyset$, indicating that $R(H)$ is transitive on $V(\Gamma_N)$. Let $u\in \Delta_i\cap H$ and $v\in \Delta_i\cap yH$. Then $R(G)$ has an element $\a$ mapping $u$ to $v$, and since $\Delta_i$ is a block, we have $\Delta_i^\a=\Delta_i$. Since $u\in H$ and $v\in yH$, $\a$ interchanges $H$ and $yH$. It follows that $|\Delta_i\cap H|=|(\Delta_i\cap H)^\a|=|\Delta_i\cap yH|$, namely $|\Delta_i\cap H|=|\Delta_i\cap yH|=\frac{1}{2}|\Delta_i|$. Thus, $R(H)_{\Delta_i}$ has two orbits on $\Delta_i$, that is, $\Delta_i\cap H$ and $\Delta_i\cap yH$, and by the semiregularity of $R(H)$, we have $|R(H)_{\Delta_i}|=\frac{1}{2}|\Delta_i|$. Since $R(H)$ is abelian and transitive on $V(\Gamma_N)$, $R(H)_{\Delta_i}$ fixes every vertex in $V(\Gamma_N)$, yielding $R(H)_{\Delta_i}\leq N$. Since $N$ is regular on $\Delta_i$, we have $|N:R(H)_{\Delta_i}|=2$.
Clearly, $R(H)_{\Delta_i}\leq R(H)\cap N$, and since $N$ fixes $\Delta_i$, we have $R(H)\cap N\leq R(H)_{\Delta_i}$. Consequently, $R(H)_{\Delta_i}=R(H)\cap N$ and $|N:R(H)\cap N|=2$.

Since $N$ is the kernel of $X$ on $V(\Gamma_N)$, we have $R(H)N/N\leq \Aut(\Gamma_N)$, and as $R(H)N/N\cong R(H)/(R(H)\cap N)$ is cyclic and acts transitively on $V(\Gamma_N)$, $R(H)N/N$ is regular on $\Gamma_N$, implying that $\Gamma_N$ is a circulant on $R(H)N/N$. Thus, $\Ga_N$ is an $(A/N,2)$-distance transitive circulant on $R(H)N/N$.

Write
$$M=\la g^2\mid g\in N\ra, \ \mbox{ and } L=\la g^2\mid g\in N\cap R(H)\ra.$$
Then $M$ is characteristic in $N$, and hence $M\unlhd X$.
Since $|N:N\cap R(H)|=2$, we have $M\leq N\cap R(H)$, and since $N\cap R(H)$ is cyclic, we have  $|N\cap R(H):L|=1$ or $2$. It follows $|N:L|=|N:N\cap R(H)||N\cap R(H):L|=2$ or $4$, and since
$L\leq M$, we have $|N:M|=2$ or $4$.
\qed

Considering the provided $\Gamma$, $X$, and $N$ according to Lemma~\ref{base}, we examine two scenarios: one where $X/N$ is quasiprimitive on $V(\Gamma_N)$ as discussed in Lemma~\ref{quasi-cover}, and the other where $X/N$ is biquasiprimitive on $V(\Gamma_N)$ as detailed in Lemma~\ref{biquasi}.

\begin{lemma}\label{quasi-cover}
Let $\Ga=\Cay(\D_{2n},S)$ be a connected $(X,2)$-distance transitive graph, where $n\geq 2$ and $R(\D_{2n})\leq X\leq\Aut(\Ga)$. Let $N$ be a non-trivial normal subgroup of $X$ such that $\Ga$ is an $N$-cover of $\Ga_N$.
If $X/N$ is quasiprimitive on $V(\Ga_N)$, then one of the following holds:
\begin{enumerate}[\rm (1)]
  \item $\Ga$ is $(X,2)$-arc transitive;
  \item  $N\cong\mz_p$ for a prime $p$, and $\Ga\cong \K_{m[b]}$ for integers $m\geq3$ and $b\geq2$ with $mb=2n$;
  \item $N\cong\ZZ_2$, $\Ga\cong\K_{p,p}-p\K_2$ and $X/N\cong\AGL(1,p)$ for a prime $p\geq5$.
\end{enumerate}
\end{lemma}

\proof Assume that $\Ga$ is not $(X,2)$-arc transitive. To finish the proof, we only need to show that part (2) or (3) of Lemma~\ref{quasi-cover} holds.
By Proposition~\ref{kernel}, $N$ is the kernel of $X$ acting on $V(\Ga_N)$, $X/N\leq \Aut(\Ga_N)$, and $\Ga_N$ is $(X/N,2)$-distance transitive but not $(X/N,2)$-arc transitive.
In particular, both $\Ga$ and $\Ga_N$ have valency at least $3$.
It is worth noting that a $2$-distance transitive group of automorphisms of a graph with  girth at least $5$ is $2$-arc transitive. Consequently, both $\Gamma$ and $\Gamma_N$ have girth less than $5$.

Let $G=\D_{2n}$ and let $t=|V(\Ga_N)|$.
Then $t\geq 4$ as $\Ga_N$ has valency at least $3$.
Since $\Ga$ is an $N$-cover of $\Ga_N$, Lemma~\ref{base} means that either $N<R(H)$ and $\Ga_N$ is a dihedrant on $R(G)/N$, or $|N:N\cap R(H)|=2$  and $\Ga_N$ is a circulant on $R(H)N/N$.
Moreover, for the latter case, there is a normal subgroup $M$ of $X$ such that $|N:M|=2$ or $4$.
Let $T=\soc(X/N)$.
We shall proceed with the proof by examining these two cases.

\medskip
{\noindent\bf Case 1.} $N<R(H)$ and $\Ga_N$ is a dihedrant on $R(G)/N$.
\medskip

In this case, we aim to prove that part~(2) of the lemma holds, that is, $N\cong\mz_p$ for a prime $p$ and $\Ga\cong \K_{m[b]}$ for integers $m\geq3$ and $b\geq2$ with $mb=2n$.

Note that $X/N$ is quasiprimitive on $V(\Ga_N)$, and contains a regular dihedral subgroup $R(G)/N$. By Proposition~\ref{Qd-group}, $X/N$ is $2$-transitive, forcing $\Ga_N\cong\K_{t}$ with $t\geq 4$. Then $|V(\Ga)|=t|N|$ and $R(G)/N\cong \D_{t}$.
Since $\Ga_N$ is $X/N$-arc transitive but not $(X/N,2)$-arc transitive, $X/N$ is $2$-transitive but not $3$-transitive on $V(\Gamma_N)$.
By Proposition~\ref{Qd-group} and Table~\ref{d-group}, one may read out $(t,X/N)$ as one of the following lists.
\begin{itemize}
   \item[List 1:] $t=4$ and $X/N\cong\A_4$;
   \item[List 2:] $t=16$ and  $X/N\cong\ZZ_2^4:\Sy_6$, $\ZZ_2^4:\A_6$, $\ZZ_2^4:\Sy_5$ or $\ZZ_2^4:\GammaL(2,4)$;
   \item[List 3:] $t=q+1$ and $X/N=\PSL(2,q).o$, where $q=r^f\equiv3\pmod{4}$ with $r$ a prime, and $o\leq\ZZ_2\times\ZZ_f$ does not contain the diagonal automorphism of $\PSL(2,q)$.
\end{itemize}

Let us first prove a claim.

\medskip
\noindent{\bf Claim:} If $N\cong \ZZ_p$ then $\Ga\cong\K_{m[b]}$ for some integers $m\geq3$ and $b\geq2$ with $mb=2n$.

\medskip
For List~1, $t=4$, $X/N\cong\A_4$, and $\Ga_N=\K_4$. Since $\Ga$ is a cyclic arc transitive cover of $\K_4$, by \cite[Theorem 6.1]{FK} we obtain that $\Ga$ is the $3$-cube $\Q_3$ or the generalized Petersen graph $P(8,3)$.
As indicated in \cite[Page 16]{CD}, $P(8,3)$ has girth $6$, and hence $\Ga\cong\Q_3$ with $\Aut(\Ga)\cong\ZZ_2^3:\Sy_3$.
Since $\Ga$ is $X$-arc transitive, we have $24\bigm| |X|$.
It is easy to see that $X=\ZZ_2^3:\Sy_3$, $\A_4\times\ZZ_2$, or $\Sy_4$ (two conjugacy classes), which can also be checked by Magma~\cite{Magma}.
Note that $N$ is a normal subgroup of $X$ of order $p$.
Then $X\not=\Sy_4$, and $X\not=\A_4\times\ZZ_2$ because $\A_4\times\ZZ_2$ has no subgroup isomorphic to $\D_8$.
Thus, $X=\ZZ_2^3:\Sy_3$ and so $\Ga$ is $(X,2)$-arc transitive, a contradiction.

For List~2, we have $T=\soc(X/N)\cong\ZZ_2^4$ and $\Ga_N\cong\K_{16}$.
Let $N\leq Y\unlhd X$ with $Y/N=T$.
Then $Y=N.T$ is transitive on $V(\Ga)$.
First assume $p=2$. Then $\Ga$ is arc transitive $\ZZ_2$-cover of $\Ga_N\cong\K_{16}$, and hence $\Ga$ has order $32$ and valency $15$. From~\cite{Conder}, there are exactly $3$ such graphs, and by Magma~\cite{Magma}, two of them are not dihedrants, and the other one is isomorphic to $\K_{16,16}-16\K_2$. Furthermore,
each arc transitive subgroup of $\Aut(\K_{16,16}-16\K_2)$ containing $Y$ as a normal regular subgroup  has no regular subgroup isomorphic to $\D_{32}$, a contradiction.

Now assume that $p$ is odd. Let $Y_2$ be a Sylow $2$-subgroup of $Y$. Then $Y_2\cong\ZZ_2^4$, and $Y=N: Y_2\cong \ZZ_p.\ZZ_2^4$. Write $C=C_Y(N)$. Then $C=N\times C_2$, where $C_2$ is a Sylow $2$-subgroup of $C$.
Applying the $N/C$-theorem (see~\cite[P.20, 4.5 Satz]{Hupp}), $Y/C$ is isomorphic to a subgroup of $\Aut(P)\cong\ZZ_{p-1}$, which implies that $Y/C$ is cyclic. As $Y_2\cong\ZZ_2^4$, we have $8\mid |C|$ and hence $C_2\cong\ZZ_2^3$ or $\ZZ_2^4$. Since $Y\unlhd X$, we have $C\unlhd X$, and since $C_2$ is characteristic in $C$, we have $C_2\unlhd X$.
Clearly, $C_2$ has at least three orbits on $V(\Ga)$. If $\Ga$ is a cover of $\Ga_{C_2}$, by Lemma~\ref{base}, we would have either $C_2\leq R(H)$ or $|C_2:C_2\cap R(H)|=2$, both of which are impossible because $R(H)$ is cyclic and $C_2\cong\ZZ_2^3$ or $\ZZ_2^4$. Thus, $\Ga$ cannot be a cover of $\Ga_{C_2}$, and since $Y$ is vertex transitive on $V(\Ga)$, Proposition~\ref{redu} implies that $\Ga\cong\K_{m[b]}$ for some integers $m\geq3$ and $b\geq2$ with $mb=2n$.

For List 3, we have that $X/N=\PSL(2,q).o$ contains the subgroup $R(G)/N\cong\D_{q+1}$.
By Lemma~\ref{proper-PSL}~(2), $R(G)/N\leq \PSL(2,q)$, and so $G\cong\D_{p(q+1)}$ as $N\cong\ZZ_p$.
Let $T=\soc(X/N)$.
Since $X/N=T.o$, there is a subgroup $Y\unlhd X$ such that $Y/N=T$.
Since $N\cong\ZZ_p$ with $p$ a prime, Proposition~\ref{mult} means that $Y=N.T$ is a central extension, $Y=NY'$, and $Y'\cong M.T$ with $M\leq\Mult(T)\cap N$.

Since $q\equiv3\pmod{4}$, it follows from \cite[P. 302, Table 4.1]{Gore-82} that $\Mult(\PSL(2,q))\cong\ZZ_2$.
If $p=|N|=2$, then \cite[Section 3.3.6]{Wilson} implies that $Y=\ZZ_2\times\PSL(2,q)$ or $\SL(2,q)$.
Note that $\D_{2(q+1)}\cong R(G)\leq Y$ as $R(G)/N\leq T=Y/N$.
Then $Y$ has an element of order $q+1$, say $z$.
By \cite[Lemma 2.9]{Pan}, $\SL(2,q)$ has no subgroup isomorphic to $\D_{2(q+1)}$, and we may let $Y=\ZZ_2\times\PSL(2,q)$.
Then $z=z_1z_2$ with $z_1\in\ZZ_2$ and $z_2\in\PSL(2,q)$, and since $q+1$ is even, $z_2$ has order $q+1$, which is impossible by~\cite[Table 8.1]{BHR}. Thus,
$p$ is an odd prime, and so $\Mult(T)\cap N=1$, compelling $M=1$.
Then $Y=N\times Y'$, and hence $N\leq Z(Y)$, the center of $Y$. This is impossible because $N<R(G)\leq Y$ and $R(G)\cong \D_{p(q+1)}$. This completes the proof of Claim.

\medskip

Now let $|N|=\ell p$ for a prime $p$. In this setting, $N$ possesses a characteristic subgroup $K$ of order $\ell$, with $|N:K| = p$. Since $N\unlhd X$, we have $K\unlhd X$. By Lemma~\ref{zp-cover}, $\Ga_K$ is an $N/K$-cover of $\Ga_N$ and $(X/K,2)$-distance transitive. Note that $N/K\cong\mz_p$. Since $\Ga$ is not $(X,2)$-arc transitive, Proposition~\ref{kernel} implies that $\Ga_K$ is not $(N/K,2)$-arc transitive. Furthermore, the quasiprimitivity of $X/N$ on $\Ga_N$ implies the quasiprimitivity of $(X/K)/(N/K)$ (identifying with $X/N$) on $(\Ga_K)_{N/K}$ (identifying with $\Ga_N$). By Claim, $\Ga_K\cong\K_{m[b]}$ for some integers $m\geq3$ and $b\geq2$ with $mb=2n$. By Proposition~\ref{Kmb}, $K=1$ because $\Ga$ is a $(K,2)$-distance cover of $\Ga_K$, and hence $\Ga=\Ga_K\cong\K_{m[b]}$ for some integers $m\geq3$ and $b\geq2$ with $mb=2n$, as required. This completes the proof of Case~1.

\medskip
{\noindent\bf Case 2.} $|N:N\cap R(H)|=2$ and $\Ga_N$ is a circulant on $R(H)N/N$.
\medskip

In this case, we aim to prove that part~(3) of the lemma holds, that is, $N\cong\ZZ_2$, $\Ga\cong\K_{p,p}-p\K_2$ and $X/N\cong\AGL(1,p)$ for a prime $p\geq5$.

Suppose that $X/N$ is not $2$-transitive on $V(\Ga_N)$. By Proposition~\ref{Qc-group}, $|V(\Ga_N)|=p$ for a prime $p$ and $X/N\leq \AGL(1,p)$. Since $\Ga_N$ is a circulant with valency at least $3$, it follows from Proposition~\ref{cir} that $\Ga_N=\K_p$ or the Paley graph $\P(p)$.
If $\Ga_N=\P(p)$, then $\Ga$ is $(X,2)$-distance transitive $N$-cover of $\P(p)$ with $|N|\geq 2$, which is impossible by Lemma~\ref{Paley}. Thus, $\Ga_N=\K_p$, and since $\Ga_N$ is $X/N$-arc transitive, $X/N=\AGL(1,p)$ is $2$-transitive on $V(\Ga_N)$, a contradiction.

Thus $X/N$ is $2$-transitive on $V(\Ga_N)$, yielding $\Ga_N\cong\K_t$.
If $X/N$ is $3$-transitive, then $\Ga_N$ is $(X/N,2)$-arc transitive, a contradiction.
Thus, $X/N$ is not $3$-transitive on $\Ga_N$. By Proposition~\ref{Qc-group}, $t$ and $X/N$ are listed in the following:

\begin{itemize}
  \item[List A:] $t=11$ and $X/N\cong\PSL(2,11)$;
  \item[List B:] $t=(q^d-1)/(q-1)$ and $X/N\cong\PGL(d,q).o$, where $d\geq3$, $q$ is a prime power and $o\leq\PGammaL(d,q)/\PGL(d,q)$;
 \item[List C:] $t=p$ and  $X/N\cong\AGL(1,p)$.
\end{itemize}

Since $R(H)N/N$ acts regularly on $V(\Ga_N)$, we have $t=|R(H)N/N|\geq 4$. Moreover, as $|R(H)N/N|\leq|R(G)N/N|=|R(G)/(R(G)\cap N)|$ and $R(G)$ is dihedral, we conclude that $N\cap R(G)<R(H)$, leading to $N\cap R(H)=N\cap R(G)$. Consequently, $R(H)N/N\cong R(H)/(R(H)\cap N)\leq R(G)/(R(H)\cap N)\cong R(G)N/N$, implying $\D_{2t}\cong R(G)N/N\leq X/N$.

Since $\PSL(2,11)$ has no subgroup isomorphic to $\D_{22}$, List~A cannot occur. Suppose that List B is applicable. Then $X/N\leq\PGammaL(d,q)$ and $R(H)N/N$ is the Singer cycle of order $(q^d-1)/(q-1)$ (see~\cite[Corollary 1.2]{Li03}). Since $|R(G)N/N:R(H)N/N|=2$, Lemma~\ref{proper-PSL}~(1) means that $R(G)N/N$ is not a dihedral group, yielding a contradiction.

Thus, we have List C, that is, $t=p$ and $X/N\cong\AGL(1,p)$. Since $p$ is prime and $p\geq 4$, we have $p\geq 5$. By Lemma~\ref{base}, $X$ has a normal subgroup $M$ such that $|N:M|=2$ or $4$.
By Lemma~\ref{zp-cover}, $\Gamma$ is a $M$-cover of $\Gamma_M$ and $\Gamma_M$ is a $N/M$-cover of $(\Gamma_M)_{N/M}$, where $(\Gamma_M)_{N/M}=\Gamma_N$ by identifying the orbit $u^N$ of $u$ under $N$ with the orbit $(u^M)^{N/M}$ of $u^M$ under $N/M$. Furthermore, $|V(\Gamma_M)|=|N/M||\Gamma_N|=p|N/M|$, and $\Gamma_M$ is $(X/M,2)$-distance transitive of valency $p-1$.

Observing that $(X/M)/(N/M)\cong X/N\cong \AGL(1,p)$ where $p\geq 5$, we find that $|X/M|=p(p-1)|N/M|$.
Let $P$ be a Sylow $p$-subgroup of $X$. Then $PM/M\cong \ZZ_p$ is a Sylow $p$-subgroup of $X/M$ and $N/M\rtimes PM/M\unlhd X/M$. Since $p\geq 5$ and $|N/M|=2$ or $4$, by the Sylow Theorem we have that $PM/M$ is normal in  $N/M\rtimes PM/M$ and hence characteristic. It follows that $PM/M\unlhd X/M$ as $N/M\rtimes PM/M\unlhd X/M$.

Suppose $|N:M|=4$. Since $\Gamma_M$ is $(X/M,2)$-distance transitive and $|V(\Gamma_M)|=p|N/M|=4p$, the quotient graph $(\Gamma_M)_{PM/M}$ has order $4$.
By Proposition~\ref{redu}, either $\Ga_M\cong\K_{m[b]}$ with $m\geq 3$ and $b\geq 2$, or $\Ga_M$ is a cover of $(\Ga_M)_{MP/M}$. The former yields constraints $mb = 4p$ and $(m-1)b = p-1$, which leads to an impossible solution. For the latter, $(\Ga_M)_{PM/M}$ has the same valency $p-1\geq 4$ as $\Ga_M$, which is impossible because $|V((\Ga_M)_{PM/M})|=4$.

Now we have $|N:M|=2$. Then $|V(\Gamma_M)|=p|N/M|=2p$, and since $PM/M\unlhd X/M$, $\Gamma_M$ is a bipartite graph of order $2p$ of valency $p-1$, forcing $\Gamma_M\cong \K_{p,p}-p\K_2$. Since $(X/M)/(N/M)\cong X/N\cong \AGL(1,p)$, $X/M$ is solvable. If $M\not=1$, then $\Ga$ is a $M$-cover of $\Ga_M$, contradicting to Lemma~\ref{2pcover}. Then $M=1$, $N\cong \ZZ_2$ and  $\Ga\cong\K_{p,p}-p\K_2$ with $X/N\cong\AGL(1,p)$ for a prime $p\geq5$, as required. This concludes the proof.\qed

Now we deal with the case when $X/N$ is biquasiprimitive on $V(\Ga_N)$.

\begin{lemma}\label{biquasi}
Let $\Ga=\Cay(\D_{2n},S)$ be a connected $(X,2)$-distance transitive graph, and let $N\unlhd X$ be such that $\Ga$ is an $N$-cover of $\Ga_N$, where $n\geq2$ and $R(\D_{2n})\leq X\leq\Aut(\Ga)$.
If $X/N$ is biquasiprimitive on $V(\Ga_N)$, then either $\Ga$ is $(X,2)$-arc transitive, or $\Ga\cong\K_{n,n}-n\K_2$ for some $n\geq 4$.
\end{lemma}

\proof Assume that  $\Ga$ is not $(X,2)$-arc transitive. To finish the proof, we will show that $\Ga\cong\K_{n,n}-n\K_2$ for some $n\geq 4$.
Let $G=\D_{2n}$ with $n\geq 2$.

By Proposition~\ref{kernel},  we ascertain that $N$ serves as the kernel of $X$ acting on $V(\Ga_N)$ and $X/N\leq \Aut(\Ga_N)$. Furthermore, $\Ga_N$ is $(X/N,2)$-distance transitive but not $(X/N,2)$-arc transitive. Since $X/N$ is biquasiprimitive on $V(\Ga_N)$, every non-trivial normal subgroup of $X/N$ has one or two orbits on $V(\Ga_N)$ and there is $M/N\unlhd X/N$ such that $M/N$ has two orbits on $V(\Ga_N)$, say $\Delta_1$ and $\Delta_2$. Set $\Delta=\{\Delta_1,\Delta_2\}$. Since $\Ga_N$ is $X/N$-arc transitive, $\Ga_N$ is a connected bipartite graph with $\Delta_1$ and $\Delta_2$ as partite sets, and since $\Ga_N$ is not $(X/N,2)$-arc transitive, $\Ga_N$ has valency at least $3$. Thus, $|\Delta_1|=|\Delta_2|\geq 3$.
Additionally, if $\Ga$ has girth at least $5$, then $\Ga$ is $(X,2)$-arc transitive, a contradiction.
Thus, $\Ga$ has girth $4$.

Let $|\Delta_1|=|\Delta_2|=t$. Since $\Ga_N$ is not $(X/N,2)$-arc transitive, by Lemma~\ref{bipart} we may assume that $\Ga_N\not\cong \K_{t,t}$. Specifically, $t=|\Delta_1|=|\Delta_2|\geq 4$.

Write $B=X/N$ and let $B^*$ be the kernel of $B$ acting on $\Delta$, that is, the subgroup of $B$ fixing $\Delta_1$ and $\Delta_2$ setwise. Then $|B:B^*|=2$.
By Proposition~\ref{biquasi-d}, $B$ is biprimitive on $V(\Ga_N)$, and hence $B^*$ has primitive actions on $\Delta_1$ and $\Delta_2$.

Let $K$ be the kernel of $B^*$ on $\Delta_1$. Since $B^*$ is primitive on $\Delta_2$, either $K$ is transitive or fixes every vertex in $\Delta_2$. For the former, $\Ga_N\cong \K_{t,t}$, a contradiction. Then the latter occurs, that is $K=1$, and so $B^*$ is faithful on $\Delta_1$.
Similarly, $B^*$ is also faithful on $\Delta_2$, and hence $B^*$ is primitive permutation groups on both $\Delta_1$ and $\Delta_2$.

Since $R(G)\leq X$ and $N\unlhd X$, Lemma~\ref{base} implies that either $N<R(H)$ and $\Ga_N$ is a dihedrant on $R(G)/N$, or $|N:N\cap R(H)|=2$ and $\Ga_N$ is a circulant on $R(H)N/N$. It follows that $B^*$ contains a regular cyclic or dihedral permutation group on $\Delta_1$, and by Propositions~\ref{Qc-group} and \ref{Qd-group}, we have: $$\mbox{ Either }B^*\leq \AGL(1,p) \mbox{, or } B^* \mbox{is $2$-transitive on }\Delta_1 \mbox{ as listed in Table~\ref{c-group} or Table~\ref{d-group}.}$$

Let $B^*\leq \AGL(1,p)$. Then $B_u=B^*_u$ for $u\in V(\Gamma_N)$ is regular on the neighbourhood $\Gamma_N(u)$. Since $|B:B^*|=2$, we deduce that $B$ is solvable and $\Ga_N$ is a $(B,2)$-distance bipartite transitive graph of order $2p$. Since $\Ga_N\not\cong\K_{p,p}$ and $\Ga_N$ has valency at least $3$, Lemma~\ref{bipartitegraphs} implies that $\Ga_N\cong\K_{p,p}-p\K_2$. By Lemma~\ref{2pcover}, $N=1$ and $\Ga=\K_{p,p}-p\K_2$, as required.

Let $B^*$ be a $2$-transitive permutation on $\Delta_1$ as listed in Table~\ref{c-group} or Table~\ref{d-group}. Let $\{u,v\}\in E(\Gamma_N)$ with $u\in \Delta_1$ and $v\in \Delta_2$.

Since $\Gamma_N$ is $B$-arc transitive, there is $\tau\in B$ such that $\tau$ interchanges $\Delta_1$ and $\Delta_2$. In particular, $\tau^2\in B^*$ and $\tau$ induces an automorphism of $B^*$. If the automorphism is an inner automorphism, that is,
there is $b\in B^*$ such that $x^\tau=x^{b^{-1}}$ for all $x\in B^*$, then  $[\tau b,B^*]=1$, and since $B=B^*\langle \tau\rangle=B^*\langle \tau b\rangle$, we have $\langle\tau b\rangle\unlhd B$. If $\tau b$ is an involution, then $\tau b$ has at least $t$ orbits with $t\geq 4$, contradicting that $B$ is biprimitive on $V(\Ga_N)$. Thus, $1\not=\langle (\tau b)^2\rangle\unlhd B^*$, and since $B^*$ is primitive on $\Delta_1$, $\langle (\tau b)^2\rangle$ has prime order, say $p$, and regular on $\Delta_1$, forcing $B^*\leq \AGL(1,p)$, which is impossible because $B^*$ is in Table~\ref{c-group} or Table~\ref{d-group}. Thus, $\tau$ induces an out automorphism of $B^*$ of order $2$.

Assume that the actions of $B^*$ on $\Delta_1$ and $\Delta_2$ are not equivalent. By Tables~\ref{c-group} or~\ref{d-group},
\begin{align*}
 (B^*,B^*_u,t)=& (\A_6,\A_5,6), (\Sy_6,\Sy_5,6), (\M_{12},\M_{11},12), (\PSL(2,11),\A_5,11), \mbox{ or } (\PGL(d,q).o,\\
 & [q^{d-1}]:\GL(d-1,q).o,(q^d-1)/(q-1)) \mbox{ with }d\geq3 \mbox{ and } q \mbox{ a prime power}.
\end{align*}

Suppose $(B^*,B^*_u,t)=(\A_6,\A_5,6)$. By the non-equivalence of actions of $B^*$, $B^*_u$ cannot fixes any vertex in $\Delta_2$, and since $B^*_u$ has an element of order $5$, $B^*_u$ is transitive on $\Delta_2$. By the Frattini argument, $B^*=B^*_uB^*_v$, and since $\{u,v\}\in E(\Gamma_N)$ and $B_u$ is transitive on $\Delta_1\backslash \{u\}$, $B_v$ is transitive on $\Delta_1\backslash \{u\}$, implying $\Gamma_N=K_{t,t}$, a contradiction. Similarly,
$(B^*,B^*_u,t)\not=(\Sy_6,\Sy_5,6)$ or $(\M_{12},\M_{11},12)$.

Suppose $(B^*,B^*_u,t)=(\PSL(2,11),\A_5,11)$. Since $\tau$ induces an out automorphism of $B^*$ of order $2$, we have $B=\PGL(2,11)$, and since $B^*_u=\A_5$, it has two orbits of length $5$ or $6$ on $\Delta_1\setminus\{u\}$. This implies that $\Ga_N\cong B(H_{11})$ or $B'(H_{11})$ (see~\cite[Example 1.1]{Pan}).
By Lemma~\ref{bipart}, $\Ga_N$ is $(B,2)$-arc transitive, a contradiction.

Suppose $(B^*,B^*_u,t)=(\PGL(d,q).o,
[q^{d-1}]:\GL(d-1,q).o,(q^d-1)/(q-1))$ with $d\geq3$ and $q$ a prime. By Proposition~\ref{Qc-group}, we may view $\Delta_1$ as the set of $1$-subspaces of $V(d,q)$, and $\Delta_2$ as the set of $(d-1)$-subspaces of $V(d,q)$ (see the paragraph after Proposition~\ref{Qc-group}). By~\cite[Example 1.2]{Pan}, $\Ga_N\cong B(\PG(d-1,q))$ or $B'(\PG(d-1,q))$, and $\Ga_N$ is $(B,2)$-arc transitive, a contradiction.

Now assume that the actions of $B^*$ on $\Delta_1$ and $\Delta_2$ are equivalent. Then $B^*_u$ fixes a vertex $w$ in $\Delta_2$, implying $B^*_u=B^*_w$.
Since $\Ga_N\ncong\K_{t,t}$, the $2$-transitivity of $B^*$ on $\Delta_1$ suggests that $\Ga_N\cong\K_{t,t}-t\K_2$, and $u$ and $w$ are not adjacent. If $B^*$ is $3$-transitive on $\Delta_1$, it is easy to see that $\Ga_N$ is $(B,2)$-arc transitive, a contradiction. Thus, $B^*$ is not $3$-transitive on $\Delta_1$ and $\Delta_2$. If $N=1$ then $\Ga=\Ga_N\cong\K_{t,t}-t\K_2$, as required. Thus, we may let $N\neq 1$.

Recall that $N<R(H)$ and $\Ga_N$ is a dihedrant on $R(G)/N$, or $|N:N\cap R(H)|=2$ and $\Ga_N$ is a circulant on $R(H)N/N$.
For convenience, write $Q=R(G)/N$ or $R(H)N/N$. Then $Q$ is cyclic or dihedral.
Let $Q^*$ be the subgroup of $Q$ fixing $\Delta_1$ and $\Delta_2$ setwise.
Then  $Q^*$ is a regular cyclic or dihedral permutation subgroup of $B^*$ on $\Delta_1$.
Write $$Q^*=\mz_t, \mbox{ or }\D_t.$$
Since $B^*$ is not $3$-transitive, one of the following holds by Tables~\ref{c-group} and~\ref{d-group}:

\begin{itemize}
  \item[(a)] $t=11$ and $(B^*,B^*_u,Q^*)=(\PSL(2,11),\A_5,\ZZ_t)$;
  \item[(b)] $t=(q^d-1)/(q-1)$ and $(B^*,B^*_u,Q^*)=(\PGL(d,q).o,[q^{d-1}]:\GL(d-1,q).o,\ZZ_t)$, where $d\geq3$ and $q$ a prime power;
  \item[(c)] $t=4$ and $(B^*,B^*_u,Q^*)=(\A_4,\ZZ_3,\D_t)$;
  \item[(d)] $t=16$, $Q^*=\D_t$ and $(B^*,B^*_u)=(\ZZ_2^4:\Sy_6,\Sy_6)$, $(\ZZ_2^4:\A_6,\A_6)$, $(\ZZ_2^4:\Sy_5,\Sy_5)$ or $(\ZZ_2^4:\GammaL(2,4),\GammaL(2,4))$;
  \item[(e)] $t=q+1$, $Q^*=\D_t$, and $(B^*,B^*_u)=(\PSL(2,q).o,\ZZ_r^f:\ZZ_{(q-1)/2}.o)$, where $q=r^f\equiv3\pmod{4}$ with $r$ a prime, and $o\leq\ZZ_2\times\ZZ_f$ does not contain the diagonal automorphism of $\PSL(2,q)$.
\end{itemize}

Let $Q$ be cyclic. Then $Q^*=\mz_t$, and cases (c)-(e) cannot occur.
For cases (a) and (b), we can regard $B$ as its inner automorphism group $\Inn(B)$ as $Z(B)=1$. Since $\tau\in B$ induces an out automorphism of $B^*$, by \cite[Lemma 2.5]{Li03} we have $C_B(Q^*)=Q^*$, which is impossible because $Q\leq C_B(Q^*)$.

Let $Q$ be dihedral. Then $Q=R(G)/N$ and $N<R(H)$. It follows that $Q^*$ is dihedral and cases (a) and (b) cannot occur.

For case (c), $\Ga_N$ has valency $3$ and $B_u=B^*_u=\mz_3$. Since $\Ga$ is a cover of $\Ga_N$,
$\Ga$ has valency $3$ and $X_\a\cong\mz_3$ for $\a\in V(\Ga)$. Thus, $X_\a$ is regular on the neighbourhood $\Ga(\a)$, and by Lemma~\ref{bipartitegraphs}, $\Ga\cong\K_{4,4}-4\K_2$, which implies that $N=1$, a contradiction.

For case (d), $Z(B^*)=1$ and $B^*$ is isomorphic to its inner automorphism group.
Since $\tau$ induces an outer automorphism of $B^*$, we may write $B^*<B\leq \AGL(4,2)$.
Thus, $B^*\not=\ZZ_2^4:\Sy_6$ or $\ZZ_2^4:\Sy_5$, because according to calculations with Magma~\cite{Magma}, $\ZZ_2^4:\Sy_6$ is a maximal subgroup of $\AGL(4,2)$ and $\AGL(4,2)$ has no subgroups of order $2|\ZZ_2^4:\Sy_5|$. Also Magma~\cite{Magma} shows that $\Out(\AGammaL(2,4))=1$, which implies that $B^*\not=\ZZ_2^4:\GammaL(2,4)=\AGammaL(2,4)$. Then $B^*=\ZZ_2^4:\A_6$ and $B=\ZZ_2^4:\Sy_6$.
Again by Magma~\cite{Magma}, $B$ does not contain a subgroup isomorphic to $\D_{32}$, a contradiction.

For case (e), $B^*=\PSL(2,q).o$ and $Q=R(G)/N\cong\D_{2(q+1)}$, where $q=r^f\equiv3\pmod{4}$ with $r$ a prime. Then $f$ is odd, and $\mz_f\times\mz_2\cong\PGammaL(2,q)/\PSL(2,q)$, which has the unique normal subgroup $\PGL(2,q)/\PSL(2,q)$ of order $2$. Since $\tau$ induces an outer automorphism of $B^*$ of order $2$, we have $\PGL(2,q)/\PSL(2,q)\leq B/\PSL(2,q)\leq \PGammaL(2,q)/\PSL(2,q)$. Thus,  $Q\leq \PGL(2,q)$ and $B=\PGL(2,q).o$.

Write $T=\PSL(2,q)$. Then $Q^*\leq T$ and $Q\not\leq T$. Let $Y\lhd X$ be such that $Y/N=T$. Since $1\not=N\leq R(H)$ is cyclic, every subgroup of $N$ is normal in $X$. By Lemma~\ref{zp-cover}, we may assume that $N\cong\ZZ_p$ for a prime $p$. By Proposition~\ref{mult}, $Y=N.T$ is a central extension, $Y=NY'$ and $Y'\cong M.T$ with $M\leq\Mult(T)\cap N$. Since $Q\not\leq T$, we have $R(G)\not\leq Y$, and since $Q^*\leq T$ and $|Q:Q^*|=2$, we have $R(G)Y/Y\cong \mz_2$.
Then $1\neq R(G)/(R(G)\cap Y)\cong R(G)Y/Y\cong\ZZ_2$, and since $R(G)\cong\D_{2p(q+1)}$, we have $R(G)\cap Y\cong\ZZ_{p(q+1)}$ or $\D_{p(q+1)}$.
Since $q\equiv3\pmod{4}$, \cite[P. 302, Table 4.1]{Gore-82}  means that $\Mult(T)\cong\ZZ_2$.
If $p=2$, then $Y\cong\ZZ_2\times\PSL(2,q)$ or $\SL(2,q)$.
By~\cite[Lemma 2.9]{Pan}, $\SL(2,q)$ has no subgroup isomorphic to $\ZZ_{2(q+1)}$ or $\D_{2(q+1)}$.
Thus $Y\cong\ZZ_2\times\PSL(2,q)$, and so $Y$ has an element of order $q+1$, say $z$.
Then $z=z_1z_2$ with $z_1\in\ZZ_2$ and $z_2\in\PSL(2,q)$, and since $q+1$ is even, $z_2$ has order $q+1$, which is impossible by~\cite[Table 8.1]{BHR}.
We therefore conclude that $p$ is an odd prime, forcing that $\Mult(T)\cap N=1$, that is $M=1$.
Then $Y=N\times T$ and $N\leq Z(Y)$, the center of $Y$.
Additionally, by~\cite[Table 8.1]{BHR}, $T$ has no element of order $q+1$, and so $N\times T$ has no subgroup isomorphic to $\ZZ_{p(q+1)}$.
This implies that $R(G)\cap Y\cong\D_{p(q+1)}$. Since $N\leq R(G)\cap Y$ and $N\leq Z(Y)$, we have $N\leq Z(R(G)\cap Y)$, which is clearly impossible.
This completes the proof.\qed

Now we are ready to demonstrate the proof of Theorem~\ref{Thm-1}.

\medskip
{\noindent\bf Proof of Theorem~\ref{Thm-1}.}
Certainly, a connected $2$-arc transitive graph is also $2$-distance transitive. Let $m\geq3$ and $b\geq2$ be integers such that $mb=2n$. It can be straightforwardly demonstrated that $\K_{m[b]}$ is a $2$-distance transitive but not $2$-arc transitive dihedrant (see~\cite[P. 186]{MP}). Thus, the sufficiency of Theorem~\ref{Thm-1} is established.

It is important to note that all $2$-arc transitive dihedrants were classified in~\cite[Theorem 1.2]{DMM}.
To establish the necessity of Theorem~\ref{Thm-1}, let $\Ga=\Cay(G,S)$ be a connected $2$-distance transitive but not $2$-arc transitive Cayley graph, where $G=\la x,y\mid x^n=y^2=1,x^y=x^{-1}\ra\cong\D_{2n}$ for some integer $n\geq2$. To complete the proof, it is sufficient to demonstrate that $\Ga\cong\K_{m[b]}$ for some integers $m\geq3$ and $b\geq2$ satisfying $mb=2n$.

Let $A=\Aut(\Ga)$.
If $A$ is quasiprimitive on $V(\Ga)$, then according to Proposition~\ref{Qd-group}, $A$ must be $2$-transitive on $V(\Ga)$. This implies $\Ga\cong\K_{2n}$, which contradicts the fact that $\Ga$ is not $2$-arc transitive.
If $A$ is biquasiprimitive on $V(\Ga)$, by Lemma~\ref{biquasi}, it would imply that $\Ga$ is $2$-arc transitive, leading to a contradiction.
Therefore, we assume that $A$ is neither quasiprimitive nor biquasiprimitive. This implies that $A$ possesses a nontrivial maximal normal subgroup $N$ with at least three orbits, and any normal subgroup of $A$ containing $N$ properly has either one or two orbits on $V(\Ga)$. Consequently, $A/N$ is either quasiprimitive or biquasiprimitive on $V(\Ga_N)$.

If $A/N$ is biquasiprimitive on $V(\Ga_N)$, then Lemma~\ref{biquasi} implies that $\Ga\cong \K_{n,n}-n\K_2$ for some $n\geq4$, or $\Ga$ is $(A,2)$-arc transitive. In both scenarios, $\Ga$ would be $2$-arc transitive, leading to a contradiction. Hence, we can assume that $A/N$ is quasiprimitive on $V(\Ga_N)$.
By Proposition~\ref{redu}, either $\Ga\cong\K_{m[b]}$ with $m\geq3$ and $b\geq 2$, or $N$ acts semiregularly on $V(\Ga)$ and $\Ga$ is a cover of $\Ga_N$.
For the former, we are done, and so we may assume that the latter occurs, that is, $\Ga$ is a cover of $\Ga_N$. By Proposition~\ref{kernel}, $A/N\leq \Aut(\Ga_N)$, $\Ga_N$ is $(A/N,2)$-distance transitive, but not $(A/N,2)$-arc transitive.
By Lemma~\ref{quasi-cover}, we reach the conclusion that either $\Ga$ is $(A,2)$-arc transitive, or $\Ga\cong\K_{m[b]}$ with $m\geq3$ and $b\geq 2$, or $\Ga\cong\K_{p,p}-p\K_2$ for a prime $p\geq5$.
It follows that $\Ga\cong\K_{m[b]}$ since $\Ga$ is not $2$-arc transitive, as desired. \qed


\begin{thebibliography}{99}

\bibitem{ACX}
B. Alspach, M.D.E. Conder, D. Maru\v{s}i\v{c}, M.Y. Xu,  A classification of $2$-arc-transitive circulants, J. Algebr. Comb.  5 (1996) 83--86.
\url{https://doi.org/10.1023/A:1022456615990}.

\bibitem{Biggs}
N.L. Biggs, Algebraic Graph Theory, Cambridge University Press, New York, 1974.
\url{https://doi.org/10.1017/CBO9780511608704}.


\bibitem{Magma}
W. Bosma, J. Cannon, C. Playoust, The MAGMA algebra system I: The user language, J. Symbolic Comput. 24 (1997) 235--265.
\url{https://doi.org/10.1006/jsco.1996.0125}.

\bibitem{BHR}
J.N. Bray, D.F. Holt, C.M. Roney-Dougal, The Maximal Subgroups of the Low-Dimensional Finite Classical Groups, Cambridge Univ. Press, 2013.
\url{https://doi.org/10.1017/CBO9781139192576}.

\bibitem{BCN}
A.E. Brouwer, A.M. Cohen, A. Neumaier, Distance-Regular Graphs, Springer Berlin, Heidelberg, 1989.
\url{https://doi.org/10.1007/978-3-642-74341-2}.


\bibitem{CJL}
J.Y. Chen, W. Jin, C.H. Li, On 2-distance-transitive circulants, J. Algebr. Comb. 49 (2019) 179--191.
\url{https://doi.org/10.1007/s10801-018-0825-3}.

\bibitem{CO}
Y. Cheng, J. Oxley, On weakly symmetric graphs of order twice a prime, J. Combin. Theory Ser. B 42 (1987) 196--211.
\url{https://doi.org/10.1016/0095-8956(87)90040-2}.


\bibitem{CJS}
B.P. Corr, W. Jin, C. Schneider, Finite $2$-distance transitive graphs, J. Graph Theory 86 (2017) 78--91.
\url{https://doi.org/10.1002/jgt.22112}.


\bibitem{Conder}
M.D.E. Conder, A complete list of all connected symmetric graphs of order $2$ to $47$,
\url{https://www.math.auckland.ac.nz/~conder/symmetricgraphs-orderupto47-byedges.txt}.

\bibitem{CD}
M.D.E. Conder, P. Dobcs\'{a}nyi, Trivalent symmetric graphs on up to $768$ vertices, J. Combin. Math. Combin. Comput. 40 (2002) 41--63.


\bibitem{Atlas}
J.H. Conway, R.T. Curtis, S.P. Norton, R.A. Parker, R.A. Wilson, Atlas of Finite Groups: Maximal Subgroups and Ordinary Characters for Simple Groups, Oxford University Press, Eynsham, 1985.


\bibitem{Dixon}
J.D. Dixon, B. Mortimer, Permutation Groups, Graduate Texts in Mathematics 163, Springer-Verlag, New York, 1996.
\url{https://doi.org/10.1007/978-1-4612-0731-3}.

\bibitem{DGLP12}
A. Devillers, M. Giudici, C.H. Li, C.E. Praeger, Locally $s$-distance transitive graphs, J. Graph Theory 69 (2012) 176--197.
\url{https://doi.org/10.1002/jgt.20574}.

\bibitem{DMM}
S.F. Du, A. Malni\v{c}, D. Maru\v{s}i\v{c}, Classification of $2$-arc-transitive dihedrants, J. Combin. Theory Ser. B 98 (2008) 1349--1372.
\url{https://doi.org/10.1016/j.jctb.2008.02.007}.


\bibitem{FK}
Y.-Q. Feng, J.H. Kwak, Cubic symmetric graphs of order a small number times a prime or a prime square, J. Combin. Theory Ser. B 97(2007) 627--646.
\url{https://doi.org/10.1016/j.jctb.2006.11.001}.


\bibitem{GR}
C.D. Godsil, G.F. Royle, Algebraic Graph Theory, Springer, New York, Berlin, Heidelberg, 2001.
\url{https://doi.org/10.1007/978-1-4613-0163-9}.


\bibitem{Gore-82}
D. Gorenstein, Finite Simple Groups: An Introduction to Their Classification, Plenum Press, New York, 1982.
\url{https://doi.org/10.1007/978-1-4684-8497-7}.


\bibitem{HFZ}
J.-J. Huang, Y.-Q. Feng, J.-X. Zhou, Two-distance transitive normal Cayley graphs, Ars Math. Contemp. 22 (2022) P.\#2.02.
\url{https://doi.org/10.26493/1855-3974.2593.1b7}.


\bibitem{Hupp}
B. Huppert, Endliche Gruppen I, Die Grundlehren der mathematischen Wissenschaften Band 134, Springer-Verlag, Berlin-New York, 1967.
\url{https://doi.org/10.1007/978-3-642-64981-3}.

\bibitem{JDLP}
W. Jin, A. Devillers, C.H. Li, C.E. Praeger, On geodesic transitive graphs, Discrete Math. 338 (2015) 168--173.
\url{https://doi.org/10.1016/j.disc.2014.11.005}.

\bibitem{JT}
W. Jin, L. Tan, Finite two-distance-transitive dihedrants, J. Aust. Math. Soc. 113 (2022) 386--401.
\url{https://doi.org/10.1017/S1446788721000409}.

\bibitem{JWZ}
W. Jin, C.X. Wu, J.X. Zhou, Two-distance-primitive graphs, Electronic J. Combin. 27 (2020) \#P4.53.
\url{https://doi.org/10.37236/8890}.


\bibitem{Jon}
G. Jones, Cyclic regular subgroups of primitive permutation groups, J. Group Theory 5 (2002) 403--407.
\url{https://doi.org/10.1515/jgth.2002.011}.


\bibitem{Kov}
I. Kov\'{a}cs, Arc-transitive dihedrants of odd prime-power order, Graphs Combin. 29 (2013) 569--583.
\url{https://doi.org/10.1007/s00373-012-1134-6}.

\bibitem{Li03}
C.H. Li, The finite primitive permutation groups containing an abelian regular subgroup, Proc. Lond. Math. Soc. 87 (2003) 725--747.
\url{https://doi.org/10.1112/S0024611503014266}.


\bibitem{LP}
C.H. Li, J.M. Pan, Finite $2$-arc-transitive abelian Cayley graphs, European J. Combin. 29 (2008) 148--158.
\url{https://doi.org/10.1016/j.ejc.2006.12.001}.



\bibitem{LSS}
C.H. Li, \'{A}. Seress, S.J. Song, $s$-Arc-transitive graphs and normal subgroups, J. Algebra 421 (2015) 331--348.
\url{https://doi.org/10.1016/j.jalgebra.2014.08.032}.


\bibitem{Mar03}
D. Maru\v{s}i\v{c}, On $2$-arc-transitivity of Cayley graphs, J. Combin. Theory Ser. B 87 (2003) 162--196.
\url{https://doi.org/10.1016/S0095-8956(02)00033-3}.

\bibitem{Mar06}
D. Maru\v{s}i\v{c}, Corrigendum to ``on $2$-arc-transitivity of Cayley graphs'' [J. Combin. Theory Ser. B 87 (2003) 162--196], J. Combin. Theory Ser. B 96 (2006) 761--764.
\url{https://doi.org/10.1016/j.jctb.2006.01.003}.


\bibitem{MP}
\v{S}. Miklavi\v{c}, P. Poto\v{c}nik, Distance-transitive dihedrants, Des. Codes Crypt. 41 (2006) 185--193.
\url{https://doi.org/10.1007/s10623-006-9008-7}.

\bibitem{Paley}
R.E.A.C. Paley, On orthogonal matrices, J. Math. Phys. 12 (1933), 311--320.
\url{https://doi.org/10.1002/sapm1933121311}.

\bibitem{Pan}
J.M. Pan, Locally primitive Cayley graphs of dihedral groups, Eur. J. Comb. 36 (2014) 39--52.
\url{https://doi.org/10.1016/j.ejc.2013.06.041}.

\bibitem{PLHL}
J.M. Pan, Y. Liu, Z.H. Huang, C.L. Liu, Tetravalent edge-transitive graphs of order $p^2q$, Sci. China Math. 57 (2014) 293--302.
\url{https://doi.org/10.1007/s11425-013-4708-8}.


\bibitem{PYZH}
J.M. Pan, X. Yu, H. Zhang, Z.H. Huang, Finite edge-transitive dihedrant graphs, Discrete Math. 312 (2012) 1006--1012.
\url{https://doi.org/10.1016/j.disc.2011.11.001}.

\bibitem{Pei}
W. Peisert, All self-complementary symmetric graphs, J. Algebra 240 (2001) 209--229.
\url{https://doi.org/10.1006/jabr.2000.8714}.

\bibitem{Q-Du-Koolen}
Z. Qiao, S.F. Du, J.H. Koolen, $2$-Walk-regular dihedrants from group divisible designs, Electronic J. Combin. 23 (2016), \#P2.51.
\url{https://doi.org/10.37236/5155}.

\bibitem{Praeger93}
C.E. Praeger, An O'Nan-Scott theorem for finite quasiprimitive permutation groups and an application to $2$-arc transitive graphs, J. London Math. Soc. s2-47 (1993) 227--239.
\url{https://doi.org/10.1112/jlms/s2-47.2.227}.


\bibitem{Schur}
I. Schur, Untersuchen \"{u}ber die darstellung der endlichen gruppen durch gebrochenen linearen substitutionen. J. Reine Angew. Math. 1904 (1904) 20--50.
\url{https://doi.org/10.1515/crll.1904.127.20}.


\bibitem{SLZ}
S.J. Song, C.H. Li, H. Zhang, Finite permutation groups with a regular dihedral subgroup, and edge-transitive dihedrants, J. Algebra  399 (2014) 948--959.
\url{https://doi.org/10.1016/j.jalgebra.2013.10.022}.


\bibitem{Wilson}
R.A. Wilson, The Finite Simple Groups, Graduate Texts in Mathematics, 251. Springer-Verlag
London, Ltd., London, 2009.
\url{https://doi.org/10.1007/978-1-84800-988-2}.


\end{thebibliography}
\end{document}